\documentclass[11pt,a4paper]{article}
\usepackage[T1]{fontenc}
\usepackage{iftex}
\ifPDFTeX\usepackage[utf8]{inputenc}\fi
\usepackage{lmodern,amsmath,amssymb,amsthm,amscd,mathtools}
\usepackage[margin=30mm]{geometry}
\usepackage{microtype,booktabs,enumitem}
\usepackage[hidelinks]{hyperref}
\usepackage[nameinlink,noabbrev]{cleveref}
\numberwithin{equation}{section}
\newtheorem{theorem}{Theorem}[section]
\newtheorem{proposition}[theorem]{Proposition}
\newtheorem{lemma}[theorem]{Lemma}
\newtheorem{corollary}[theorem]{Corollary}
\theoremstyle{definition}

\newcommand{\Q}{\mathbb Q}
\newcommand{\C}{\mathbb C}
\newcommand{\Z}{\mathbb Z}
\newcommand{\PP}{\mathbb P}
\newcommand{\Gm}{\mathbb G_m}
\newcommand{\one}{\mathbf1}

\newcommand{\dt}{\mathbin{\otimes^{\mathbf L}}}
\newcommand{\sh}[1]{\langle #1\rangle}
\newcommand{\U}{\mathcal U}
\newcommand{\Mg}{\mathcal M_g}
\newcommand{\Jac}{\mathcal J^d_g}
\newcommand{\Pencils}{\mathcal G^1_{d,g}}
\newcommand{\Kh}{\widehat K}
\DeclareMathOperator{\Pic}{Pic}

\DeclareMathOperator{\End}{End}

\DeclareMathOperator{\Ext}{Ext}
\DeclareMathOperator{\Tor}{Tor}
\DeclareMathOperator{\Gr}{Gr}
\DeclareMathOperator{\Conf}{Conf}
\DeclareMathOperator{\UConf}{UConf}
\DeclareMathOperator{\Hur}{Hur}
\DeclareMathOperator{\CHur}{CHur}
\DeclareMathOperator{\cofib}{cofib}

\DeclareMathOperator{\coker}{coker}
\DeclareMathOperator{\codim}{codim}

\DeclareMathOperator{\Crit}{Crit}
\DeclareMathOperator{\Endw}{end}
\DeclareMathOperator{\PGL}{PGL}

\DeclareMathOperator{\rk}{rk}
\DeclareMathOperator{\divisor}{div}
\setlist{itemsep=3pt,topsep=5pt}
\hypersetup{pdftitle={The rational Picard group of simple Hurwitz spaces},pdfauthor={Hailin Wen}}
\title{The rational Picard group of simple Hurwitz spaces}
\author{Hailin Wen}

\date{}
\begin{document}
\maketitle
\begin{abstract}
We prove that every line bundle on the moduli stack of connected,
simply branched covers of a rational curve is torsion, for every degree
$d\ge2$ and source genus $g\ge0$. The target is unparameterized and
the branch values are unordered. Our proof combines homological
stabilization with the geometry of linear series. For $d\ge6$, we
obtain a uniform range $g\ge9$ by studying complete frames of disk
covers. The intersection form on the compactified frame surface
controls the Leray differentials of an orbit configuration space;
insertion and collision of branch points then give bounds on the first
two stabilization kernels. Together with the results of
Deopurkar--Patel and Mullane, this leaves three pairs $(d,g)$.
Universal pencils settle $(6,7)$ and $(7,8)$, while a family of pencils
on a K3 surface settles $(6,8)$. In the geometric argument we keep
track of line bundles through their generic restrictions, so that the
conclusion concerns the full Picard group, including its algebraically
trivial part.
\end{abstract}
\tableofcontents
\section{Introduction}
\label{sec:introduction}

Fix integers $d\ge2$ and $g\ge0$. A point of $\U_{d,g}$ is a finite map
$f:C\to P$, where $C$ is a smooth connected projective curve of genus $g$,
$P$ is a smooth curve of genus zero, and $f$ has degree $d$ and simple
ramification with distinct branch values. Isomorphisms may act on both
$C$ and $P$. No branch value or sheet is labelled. We work over $\C$ and
write $\Pic(-)_\Q=\Pic(-)\otimes_\Z\Q$.

The rational Picard rank conjecture asks whether every line bundle on
$\U_{d,g}$ is torsion. The distinction between this assertion and a
calculation of numerical divisor classes matters: a bound on the
N\'eron--Severi group alone need not control algebraically trivial line
bundles. Our aim is to show that each line bundle becomes trivial after
taking a positive tensor power. Equivalently, we seek vanishing after
tensoring the entire Picard group with $\Q$.

\begin{theorem}\label{thm:main}
For every $d\ge2$ and $g\ge0$,
\[
                    \Pic(\U_{d,g})_\Q=0.
\]
For $2g+2d-2\ge3$, the same assertion holds for the coarse moduli space.
The case $(d,g)=(2,0)$ is treated separately in \cref{sec:completion}.
\end{theorem}

The proof uses the following established results. Deopurkar--Patel
prove the low-degree cases $d\le5$ \cite[Theorem A]{DP}.
Mullane proves vanishing when $d\ge g$ \cite[Theorem 1.1]{Mullane}.
Landesman--Levy prove the result for each fixed degree and sufficiently
large genus \cite[Theorem 1.2.3 and \S7]{LL}. Their componentwise stable
homology calculation is an input here. We make the genus bound uniform
in $d$ by estimating the kernels of stabilization.

There are two parts to the proof. The first compares a cover with the
cover obtained by adding two simple branch points. Once the genus is
large enough, the relevant homology maps are injective, so the known
stable calculation determines the low-degree cohomology. The second
part concerns the few small genera left over: there we regard a cover
as a pencil on its source curve and study line bundles on the space of
all pencils. Keeping these two viewpoints separate will also clarify
where topology enters and where an algebraic Picard-group argument is
needed.

\begin{theorem}[Uniform topological range]\label{thm:uniform}
If $d\ge6$ and $g\ge9$, then $\Pic(\U_{d,g})_\Q=0$.
\end{theorem}

\begin{theorem}[The geometric range needed here]\label{thm:geometric}
For $d\ge6$, rational Picard vanishing holds when
\[
                 d+1\le g\le2d-4.
\]
\end{theorem}

Only three instances of \cref{thm:geometric} are needed to deduce the
main theorem. Indeed, if $d\ge6$, $g\le8$ and $g>d$, then
\[
                    (d,g)\in\{(6,7),(6,8),(7,8)\}.
\]
The argument for $g\le2d-5$ is uniform and relatively short. The separate
argument for $g=2d-4$ is presented in its natural generality; its required
instance is $(6,8)$.

\begin{center}
\begin{tabular}{@{}lll@{}}
\toprule
Range & Method & Location\\
\midrule
$2\le d\le5$, $g\ge3$ & Deopurkar--Patel & \cite{DP}\\
$d\ge g$ & Mullane & \cite{Mullane}\\
$d\ge6$, $g\ge9$ & Projectors and stabilization & \S\ref{sec:hurwitz}--\S\ref{sec:bounds}\\
$(6,7),(7,8)$ & Universal pencils, $\rho\ge3$ & \S\ref{sec:pencils}\\
$(6,8)$ & K3 test family, $\rho=2$ & \S\ref{sec:k3}\\
\bottomrule
\end{tabular}
\end{center}

\subsection{The topological argument in outline}
Put $r=d-1$, and let $T_d$ be the transpositions in $S_d$. The algebra
$A=C_*(\Hur^{T_d};\Q)$ contains all framed disk covers, including
disconnected ones. Its bimodule $M=C_*(\CHur^{T_d};\Q)$ contains the
connected covers. Both retain every boundary product. The grading
$N$ counts branch points; it is independent of homological degree.

Inserting an equal pair $(e,e)$ gives a square operator $x_e$ of word
degree two. Set
\[
 z=\sum_{e\in T_d}x_e,\qquad \Delta=(x_e-x_f),\qquad J=(x_e).
\]
To measure the failure of stabilization, set
\[
 C_z=\cofib(z:A\sh2\longrightarrow A),\qquad X=M\dt_A C_z,
\]
where $\sh2$ raises word degree by two. The homology of this derived
cofiber records both cokernels and shifted kernels. These occur as the
two terms of a short exact sequence, rather than as a prescribed direct
sum. In particular $z^2H_i(X)=0$; the stronger assertion $zH_i(X)=0$
does not follow from the cofiber sequence.

The geometry behind the new annihilation bound is as follows. The complete
frames of a connected disk cover form a connected principal $S_d$-cover
$S$ of a punctured disk. Let
\[
 O_m(S)=\{(y_1,\ldots,y_m)\in S^m:G y_i\cap G y_j=\varnothing\ (i\ne j)\},
 \qquad P_m=\prod_{j=1}^m(1-e_j),
\]
where $e_j$ averages the action of $G=S_d$ in coordinate $j$.
We prove
\begin{equation}\label{eq:intro-projectors}
 \begin{aligned}
 P_2H^{\le1}(O_2(S);\Q)&=0 &&(N>2r),\\
 P_3H^{\le2}(O_3(S);\Q)&=0 &&(N\ge2r+4).
 \end{aligned}
\end{equation}
The projector $P_m$ keeps the part that is nontrivial in every frame
coordinate. Terms involving fewer than $m$ coordinates disappear
immediately. For $m=3$, two types of terms remain in total degree two,
and the proof consists of showing that their Leray differentials are
injective.

The compact completion $\overline S$ is essential. Every nontrivial
irreducible $S_d$-representation occurs in $H^1(\overline S)$ with
multiplicity at least four when $N\ge2r+4$. The compact intersection
form detects graph-diagonal classes and makes a three-position tensor
map injective. The punctured surface itself has $H^2(S)=0$ and has no
nondegenerate intersection pairing on its ordinary $H^1$.

Complete-frame insertions transfer \eqref{eq:intro-projectors} to
operators on Hurwitz homology. A local collision calculation improves
the first cofiber degree. The geodesic algebra $B=H_0(A)/J$ then supplies
derived coefficient bounds. Write
\[
 \Kh_i=\ker\bigl(z:H_i(M)\sh2\longrightarrow H_i(M)\bigr).
\]
The resulting kernel estimates are
\begin{equation}\label{eq:intro-kernels}
 \Endw\Kh_1\le2r+8,\qquad \Endw\Kh_2\le2r+18.
\end{equation}
Here $\Kh_{i,N+2}$ is the kernel with \emph{source} word degree $N$.
Since a genus-$g$ cover of the sphere has $N=2g+2r$ branch points,
$N+2>2r+18$ first holds at $g=9$.

\subsection{The geometric argument in outline}
Let $\Pencils$ be the scalar-rigidified stack of all degree-$d$ pencils
on smooth genus-$g$ curves. Base points and incomplete pencils are
included. Here $\Mg$ is the moduli stack of smooth genus-$g$ curves.
The map $\Pencils\to\Mg$ is proper, and the total stack is smooth.
Write $\rho=2d-g-2$, the expected dimension of the space of pencils
on a fixed curve.

For $\rho\ge3$, the Fulton--Lazarsfeld incidence comparison gives the
low-degree cohomology of every fibre, including special fibres.
Globalizing this comparison gives a three-generator rational Picard
group. Three line-bundle relations vanish on the simple Hurwitz
open and kill these generators.

At $\rho=2$, the same comparison gives only an injection in degree two.
We instead restrict global line bundles to a family of pencils on a
Picard-rank-one K3 surface. A sheaf-moduli calculation bounds the Picard
group of that whole family. Flatness of the fibre Chern class of a
global line bundle, followed by the strong Franchetta theorem, controls
its generic restriction. A conormal--Wronskian estimate controls the
remaining vertical divisors. This treats $\Pic^0$ as well as divisor
classes.

\subsection{Conventions and reading guide}
The chain complexes in the topological part are homologically graded,
with differential of degree $-1$. Deformation complexes and derived
direct images in the geometric part use cohomological grading.
Word shifts are
\[
                       (Y\sh a)_N=Y_{N-a}.
\]
For a word-graded vector space, $\Endw V$ is the supremum of the degrees
where it is nonzero. Tensor products over $A$ are derived unless stated
otherwise. The product convention for a Hurwitz move is
$(a,b)\mapsto(b,b^{-1}ab)$, and $a^b=b^{-1}ab$.

The topological and geometric parts can be read independently.
\Cref{sec:completion} gives the final coverage argument.

\section{Hurwitz spaces, square operators and the cofiber}
\label{sec:hurwitz}

The stabilization operation is easiest to organize for covers of a disk.
Choose a standard system of paths from the boundary to the branch
points. A numbering on a boundary fibre then records the cover by a
tuple of local monodromies, each a transposition. Changing the paths
gives the Hurwitz action. Concatenating tuples corresponds to gluing
two disks along boundary intervals. We include disconnected covers in
the algebra so that individual branch points can be added freely, and
keep connected covers in a bimodule over that algebra.

\subsection{The algebra and its connected bimodule}
Let $D$ be a closed disk and let $D^\circ$ be its interior. Write
$\UConf_N(D^\circ)$ for the space of unordered configurations of $N$
distinct points in $D^\circ$. Let $\Hur_N^{T_d}$ be its covering with fibre
$T_d^N$, with monodromy the Hurwitz action of the braid group $B_N$.
Equivalently, it is the classifying space of the action groupoid
$B_N\curvearrowright T_d^N$. Its points classify degree-$d$ simple disk
covers with a numbered fibre at a boundary basepoint. A tuple belongs
to $\CHur_N^{T_d}$ when its entries generate $S_d$.

For transpositions, transitivity of the generated subgroup is equivalent
to generation of $S_d$: the connected components of the graph with an
edge for each transposition are precisely the orbits, and each connected
component generates its full symmetric group. Thus $\CHur$ is the
space of connected disk covers.

Concatenation makes
\[
 A=C_*(\Hur^{T_d};\Q),\qquad M=C_*(\CHur^{T_d};\Q)
\]
an augmented algebra and an $A$-bimodule. One may use Moore rectangles,
with a numbered lower boundary interval, to make the product associative.
Appending any cover to a connected cover preserves connectedness.
Normalized chains of the action-groupoid model have word-zero part
$A_0=\Q$ and no negative homological or word degrees. The geometric
and normalized-chain models give equivalent derived module categories.

We retain all products $\sigma\in S_d$. Restricting to product one at
this point would destroy the module structure under single-letter
multiplication. Product one is imposed only in the final comparison
with covers of $\PP^1$.

\subsection{Components and the geodesic quotient}
Let $R=H_0(A)$, and let $\ell(\sigma)=d-c(\sigma)$ be the minimal number
of transpositions in a factorization of $\sigma$, where $c(\sigma)$
counts all cycles, including fixed points. Thus $\ell(\sigma)\le r=d-1$.

We use the transposition orbit theorem: a Hurwitz orbit is determined by
the product, the generated Young subgroup, and the number of factors
in each of its nontrivial blocks. For a connected block, the shortest
admissible factorization of a specified product has length
$2(s-1)-\ell(\sigma|_{\mathrm{block}})$; all longer admissible lengths
differ by even integers. See \cite[Theorem 2.5]{LewisWang} for the orbit
classification and \cite[Propositions 7.11 and 7.13]{Bianchi} for its
algebraic formulation. The support partition is needed here because we
have included disconnected covers: specifying the total product alone
does not record how the letters are distributed among the components.

For $e\in T_d$, put $x_e=[e,e]\in R_2$. Moving a product-one packet
past a single strand gives centrality of $x_e$ in $R$. Set
\begin{equation}\label{eq:ideals}
 z=\sum_{e\in T_d}x_e,\qquad
 \Delta=(x_e-x_f:e,f\in T_d),\qquad J=(x_e:e\in T_d).
\end{equation}
Here $\Delta$ and $J$ are two-sided ideals. Since
$q=|T_d|=\binom d2$ is invertible in $\Q$, we have $J=(z,\Delta)$.

\begin{lemma}\label{lem:ring}
Choose a shortest factorization $u_\sigma$ for each $\sigma$. If $Z\subset R$
is the subalgebra generated by the squares, then
\[
 R_\sigma=Z u_\sigma.
\]
Moreover, in the product-$\sigma$ sector of the connected module,
\begin{equation}\label{eq:hzero}
 H_0(M)_\sigma\cong\Q[t]\sh{2r-\ell(\sigma)},\qquad
 x_e=t,\quad z=qt.
\end{equation}
If $\ell(\sigma)=r$, all squares agree already on $R_\sigma$.
\end{lemma}
\begin{proof}
A shortest factorization has the cycles of $\sigma$ as its support blocks.
To obtain a prescribed larger support partition, join these cycles
inside each desired block by doubled tree edges. Further equal pairs
increase the prescribed block lengths by two. The resulting factorization
has the same product, support and block lengths as the given one, so the
orbit theorem identifies their classes. This proves the spanning claim
without asserting that $Z$ is a polynomial ring.

In the connected case there is one orbit for each length of the parity
of $\ell(\sigma)$ starting at $2r-\ell(\sigma)$, and insertion of any
equal pair moves to the next such length. This proves \eqref{eq:hzero}.
If $\ell(\sigma)=r$, then $\sigma$ is a $d$-cycle, so every factorization
of it is already connected. The last statement follows.
\end{proof}

The quotient
\begin{equation}\label{eq:geodesic-algebra}
 B=R/J=\bigoplus_{\sigma\in S_d}\Q u_\sigma,
 \qquad |u_\sigma|=\ell(\sigma),
\end{equation}
has multiplication $u_\sigma u_\pi=u_{\sigma\pi}$ when lengths add,
and zero otherwise. Indeed, a nonminimal factorization contains a square
in the spanning description; minimal factorizations with the same product
belong to a single orbit. This is the geodesic symmetric-group algebra
of \cite[Proposition 7.14]{Bianchi}. In particular,
\begin{equation}\label{eq:koszul}
 \Endw B=r,\qquad B\text{ is Koszul},\qquad
 \Tor_p^B(\Q,\Q)\text{ is concentrated in word }p.
\end{equation}

\subsection{The derived cofiber and shifted kernels}
Each square $x_e$ is central in the monoid of Hurwitz components.
By \cite[Lemma 3.2.3]{LL}, its chain-level multiplication refines to an
$A$-bimodule map. Choose these refinements and add them to obtain a
refinement of $z$. Define
\[
 C_z=\cofib(z:A\sh2\longrightarrow A),\qquad
 X=M\dt_A C_z,\qquad V_i=H_i(X).
\]
In this tensor product $M$ is used as a right $A$-module. Its left
$A$-action induces a left action on $X$, and hence an $R$-action on
each $V_i$. This is the action used in the coefficient estimates below.
Let
\[
 C_i=\coker\bigl(z:H_i(M)\sh2\to H_i(M)\bigr),\qquad
 \Kh_i=\ker\bigl(z:H_i(M)\sh2\to H_i(M)\bigr).
\]
The long exact sequence gives $R$-linear short exact sequences
\begin{equation}\label{eq:cone-short}
             0\longrightarrow C_i\longrightarrow V_i
                \longrightarrow\Kh_{i-1}\longrightarrow0.
\end{equation}
The square actions commute with this sequence. On $M$, left and right
square actions agree on homology: moving a product-one packet across the
old cover gives the required natural homotopy of insertion maps. Thus
the $z$ appearing in the cofiber agrees with the sum of the square
operators used below.

Both outer terms of \eqref{eq:cone-short} are killed by $z$, whence
\begin{equation}\label{eq:z-square}
                        z^2V_i=0.
\end{equation}
The extension need not split, and it does not follow that $zV_i=0$.
By \eqref{eq:hzero}, multiplication by $z$ on $H_0(M)$ is injective.
Consequently
\begin{equation}\label{eq:low-cone}
 \Kh_0=0,\qquad V_1=C_1,\qquad zV_1=0,\qquad
 JV_0=0,\qquad \Endw V_0\le2r.
\end{equation}

We will need two further inputs from Landesman--Levy. The augmentation
module $\Q$ is concentrated in word zero, with all positive words acting
by zero. Their estimate for the derived tensor of two augmentation
modules \cite[Lemma 3.2.8]{LL}, with the full transposition rack, gives
\begin{equation}\label{eq:augmentation}
          \Endw H_i(\Q\dt_A\Q)\le i.
\end{equation}
Their componentwise stable comparison
\cite[Theorems 1.4.1 and 1.4.2]{LL} identifies the rational homology of
connected transposition Hurwitz spaces, at sufficiently large word degree
for each fixed $d$ and homological degree, with that of unordered plane
configurations. In degrees one and two these stable groups are $\Q$ and
$0$. The comparison is compatible with equal-pair insertions; $z$ acts
as the nonzero scalar $q$ times ordinary two-point stabilization.
We use this comparison as the eventual destination of stabilization.
The explicit range in this paper will come from bounding the kernels
of the maps leading to it.

\section{Connected arcs and derived geodesic coefficients}
\label{sec:arcs}

This section proves the two coefficient estimates used in the numerical
argument. The first keeps the connected cover throughout the resolution.
The second uses the finite geodesic algebra in place of the augmentation.
The geometric point is that a tether lets us separate off one branch
point. Several tethers give the successive terms of a resolution, while
forgetting a tether gives a face map. For connected covers we must
control how much of this resolution remains after requiring the cut
cover to stay connected.

\subsection{Arcs that preserve every source component}
Let $f:F\to D$ be a simple cover of degree $d$ with $N$ branch points and
$c$ source components. Every component has boundary. Riemann--Hurwitz gives
\[
 \xi(f)=N-d+c=\sum_{F_a\subset F}(1-\chi(F_a))\ge0.
\]
Let $I$ be a nonempty finite union of disjoint open boundary intervals.
A tether is an embedded arc from $I$ to a branch point, with interior
disjoint from the branch points and the boundary. A system of tethers
has disjoint images and distinct terminal points; initial endpoints
move in their specified intervals under isotopy.

Cutting $q$ tethers removes $q$ ramified one-handles. If the resulting
source has $c+s$ components, then
\begin{equation}\label{eq:cut-complexity}
                    \xi(f_S)=\xi(f)-q+s.
\end{equation}
Call a system good if $s=0$, and write $K(f;I)$ for the complex of good
systems. A subsystem of a good system is good, since restoring handles
cannot split components.

We use the following cut description of links. After fixing a system and
cutting along it, disjoint additional tethers are tethers in the residual
disk with initial endpoints in the image of $I$ minus the cut endpoints.
This is again a finite union of intervals. The newly exposed slit sides
are not added to the allowed initial set. The correspondence respects
isotopies: an isotopy can be made disjoint from the fixed slits by
successively removing innermost unpunctured bigons. Conversely a residual
isotopy glues to an original one.

\begin{lemma}\label{lem:full-arcs}
The unrestricted tether complex with $N\ge1$ terminal branch points and
initial set $I$ is $(N-2)$-connected.
\end{lemma}
\begin{proof}
For one interval the unrestricted complex is contractible; this is the
arc complex used in \cite[\S4.1]{BM}. Fix one component $I_0\subset I$.
For a simplicial $h$-sphere map, $h\le N-2$, choose a maximal domain
simplex whose tethers all begin outside $I_0$. If it has $q_0$ vertices
and $q$ distinct image tethers, then $q\le q_0\le h+1\le N-1$.
Its link begins in $I_0$. Cutting the $q$ tethers leaves at least one
branch point and the intact interval $I_0$, so this link fills in the
contractible residual one-interval complex. Replace the star by the join
of its boundary with this filling. All new bad simplices are proper faces
of the old one. Descending dimension and then the number of maximal bad
simplices terminates. The resulting sphere lies in the $I_0$ complex.
A repeated image on the link would give a larger bad domain simplex,
so repeated images do not affect the cut argument.
\end{proof}

\begin{proposition}[Connected-complement connectivity]\label{prop:cm-arcs}
The complex $K(f;I)$ is $(\xi(f)-2)$-connected. The link of a good
$q$-tether system is $(\xi(f)-q-2)$-connected.
\end{proposition}
\begin{proof}
After cutting a system $S$, form its inverse-gluing graph: vertices are
the residual source components, and edges are the removed ramified
handles. Restoring all edges recovers the original components. Call a
nonempty system bad if every edge in this graph is a nonloop.

Four elementary properties govern the proof. A nongood system has a bad
face, obtained by restoring its loop edges. Further cuts cannot turn an
old nonloop into a loop. If additional cuts split a residual component,
their bad face together with the old bad system is a larger bad system.
Finally, if the additional cuts preserve all residual components, their
edges are loops, and remain loops after any old handles are restored.
They are statements about components of the cut surface, and so remain
valid when transporting tethers changes their monodromy labels.

Induct on $N$, allowing arbitrary $d,c,I$. The claim is vacuous for
$\xi=0$. For $\xi\ge1$, cut a full tether system to obtain $d$ disks.
The inverse-gluing graph has cycle rank $\xi$, hence has a nonbridge
edge. Removing only that handle preserves the original components and
gives a vertex of $K$.

For $\xi\ge2$, extend an $h$-sphere in $K$, $h\le\xi-2$, to an
$(h+1)$-ball in the unrestricted complex, using \cref{lem:full-arcs}
and $\xi\le N$. Choose a bad domain simplex of maximal dimension. It
does not lie in the boundary. Let it have $q_0$ vertices and $q$ distinct
image tethers; cutting gives $c+s$ components with $s\ge1$.
Maximality and the four properties imply that its link maps to the
good complex of this cut cover. The link has dimension
\[
 h-q_0+1\le h-q+1\le\xi-q-1\le\xi(f_S)-2.
\]
It therefore fills by induction. Star replacement introduces no bad
simplex containing a new vertex; the remaining bad simplices are proper
faces of the one removed. Repetition gives a filling entirely in $K$.
If an image tether occurred again on the link, adjoining that vertex
would contradict maximality, so degeneracies cause no additional case.
For a good system, $s=0$ in \eqref{eq:cut-complexity}; the link statement
follows from the same cut description.
\end{proof}

We next translate the connectivity statement into a bound for a derived
tensor product. There are two comparisons to make: the usual tether
resolution has a free right Hurwitz action, and tensoring it with
connected covers selects precisely the tethers with connected residual
cover. We describe the comparison before using its connectivity.

\begin{proposition}[Connected derived coefficient bound]\label{prop:cm}
For $r=d-1$,
\begin{equation}\label{eq:cm}
             \Endw H_i(\Q\dt_A M)\le r+i.
\end{equation}
The reflected statement holds on the other side.
\end{proposition}
\begin{proof}
Let $S_\bullet\to\Hur$ be the geometric unrestricted tether resolution
of \cite[\S\S4.2--4.3]{BM}, with its right concatenation action. The
standard straight tethers give right-module equivalences
\[
                    S_p\simeq T_d^{p+1}\times\Hur,
\]
where the generators have word degree $p+1$. In word $N$ this is the
stabilizer equivalence: braids act transitively on tether systems, and
the stabilizer is $B_{N-p-1}$ acting on the residual cover. Retain the
geometric face maps and their coherences. In degree zero there is one
tether and a residual cover with $N-1$ branch points; in degree one
there are two tethers and a residual cover with $N-2$ branch points.
The two faces forget one tether at a time, with the sheet labels
transported by the corresponding Hurwitz move.

Let $K_p$ parametrize covers with $p+1$ tethers whose residual cover is
connected. Concatenation gives a balanced map from the two-sided bar
construction
\[
                    B(S_p,\Hur,\CHur)\longrightarrow K_p.
\]
A bar simplex $(s,a_1,\ldots,a_k,m)$ maps to its concatenation, retaining
the original tethers. This commutes with bar faces and with forgetting a
tether. The free-module equivalence identifies the source with
$T_d^{p+1}\times\CHur$; the same stabilizer argument identifies its
map to $K_p$ as an equivalence. Thus the comparison is natural on the
whole augmented diagram. For each face $d_i$ this gives a commutative
comparison
\[
\begin{CD}
 B(S_p,\Hur,\CHur) @>>> K_p\\
 @V{d_i}VV @VV{d_i}V\\
 B(S_{p-1},\Hur,\CHur) @>>> K_{p-1}.
\end{CD}
\]
Here forgetting a tether may reattach its ramified handle to the
residual cover. Since reattaching a handle cannot disconnect that
cover, the right-hand face map stays inside the connected resolution.

The cofiber of the unrestricted augmentation is $\Q$: in positive words
the full arc complex is contractible, and word zero contains only the
augmentation unit. Tensoring this diagram with $M$ shows that the
cofiber of $|K_\bullet|\to\CHur$ computes $\Q\dt_A M$.
Over a connected cover in word $N$, its homotopy fibre is $K(f;I)$,
which is $(N-r-2)$-connected by \cref{prop:cm-arcs}. The augmentation
is $(N-r-1)$-connected, so its homology cofiber vanishes for
$i\le N-r-1$, equivalently for $N>r+i$. Reflection reverses
concatenation and proves the other-side assertion.
\end{proof}

\subsection{The geodesic coefficient estimate on both sides}
Passing from $R$ to $B=R/J$ removes every factorization that contains an
equal pair. Only minimal transposition factorizations remain. The
Koszul property in \eqref{eq:koszul} says that a Tor class in homological
degree $p$ over $B$ has word degree exactly $p$. We combine this precise
grading with the augmentation estimate \eqref{eq:augmentation}.

\begin{proposition}\label{prop:geo}
For $B=R/J$,
\begin{equation}\label{eq:geo}
 \Endw H_j(\Q\dt_A B)\le2j,
 \qquad \Endw H_j(B\dt_A\Q)\le2j.
\end{equation}
\end{proposition}
\begin{proof}
Tensor the unrestricted tether resolution on the right with the
discrete left module $B$. All levels are now discrete, so its homology
is computed by the normalized complex
\[
 D_j=\Q[T_d^j]\otimes B,\qquad D_0=B,
\]
with differential
\begin{equation}\label{eq:rack-differential}
 d(a_1,\ldots,a_j;b)=
 \sum_{h=1}^j(-1)^{h-1}
 (a_1,\ldots,\widehat a_h,\ldots,a_j;
                         a_h^{a_{h+1}\cdots a_j}b).
\end{equation}
The augmentation term has homological degree zero. Hence this complex
computes $U_j=H_j(\Q\dt_A B)$, with its positive \emph{right} $B$-action.
The relation $ba=ab^a$ gives a homotopy for right multiplication by $a$:
\[
 H_a(a_1,\ldots,a_j;b)=(-1)^j(a_1,\ldots,a_j,a;b^a).
\]
In $dH_a+H_ad$, deletion of the last label is multiplication by $a$;
the other terms cancel in pairs with the signs in
\eqref{eq:rack-differential}. Thus the positive right action on $U_j$
is zero. This is the explicit form of the twist homotopy of
\cite[\S5.2]{BM} in our grading convention.

The connective change-of-rings spectral sequence is
\[
 E^2_{p,q}=\Tor_p^B(U_q,\Q)
       =U_q\otimes\Tor_p^B(\Q,\Q)
       \ \Longrightarrow\ H_{p+q}(\Q\dt_A\Q).
\]
Its differential is $(p,q)\mapsto(p-s,q+s-1)$. We have $U_0=\Q$.
Assume the desired bound below $j$. There is no outgoing differential
from $(0,j)$, and incoming terms have indices $(s,j-s+1)$, with
$2\le s\le j+1$. By Koszulity their words are at most
\[
                  s+2(j-s+1)\le2j.
\]
Consequently any $U_j$ class of word greater than $2j$ survives to the
abutment, contradicting \eqref{eq:augmentation}. This proves the first
bound by induction. The first-quadrant filtration is finite in each
homological degree.

For the second bound, reverse tuples and reflect the disk. On braid
generators this is $\sigma_i\mapsto\sigma_{N-i}^{-1}$; it intertwines
the chosen Hurwitz convention and reverses concatenation. It preserves
$J$, the augmentation and word length, inducing the chain anti-equivalence
with the graded opposite algebra. On $B$ it sends $u_\sigma$ to
$u_{\sigma^{-1}}$. Thus it transports the first derived tensor
to the second. Explicitly, reversal on the normalized length-$j$ complex
has sign $(-1)^{j(j-1)/2}$; this identifies the differentials and confirms
that the action contracted on the opposite side is the positive left
action.
\end{proof}

\begin{corollary}\label{cor:b-cost}
If $W$ is a discrete word-graded left $R$-module whose word support is
bounded below and contained in degrees at most $h$, then
\[
 \Endw H_p(B\dt_A W)\le h+2p.
\]
Moreover
\begin{equation}\label{eq:bcm}
                 \Endw H_i(B\dt_A M)\le2r+i.
\end{equation}
\end{corollary}
\begin{proof}
Filter $W$ by decreasing word degree. Each quotient is an augmentation
module concentrated in a single word. The two support bounds make the
filtration finite, and
\cref{prop:geo} bounds its derived tensor pieces by $h+2p$. Extensions
do not increase support beyond the maximum of their pieces. For
\eqref{eq:bcm}, filter the discrete right module $B$ in the same way.
Its words lie between zero and $r$. Tensoring with $M$ and using
\cref{prop:cm} gives $r+(r+i)$. The resulting object retains the right
$A$-action inherited from $M$.
\end{proof}

\section{Complete frames and the three-coordinate projector}
\label{sec:projectors}

We now study the extra data used when inserting branch points: an
unbranched point together with a numbering of all sheets above it.
Several such choices give an orbit configuration space on the surface
of complete frames. The aim is to show that its low-degree cohomology
has no part on which every frame coordinate acts nontrivially.

There are two geometric ingredients. The punctured frame surface carries
the free group action needed to define orbit configurations. Its compact
completion carries the nondegenerate intersection form needed to detect
diagonal classes. We first relate these two surfaces, then use the
pairing to compute the relevant differentials.

\subsection{The compact complete-frame surface}
Fix a connected degree-$d$ disk cover with $N$ branch points and boundary
product $\sigma$. Over an unbranched interior point $t$, take the set of
all bijections $\{1,\ldots,d\}\to f^{-1}(t)$. These sets form a principal
$G=S_d$-cover $S$ of the punctured open disk. The monodromy is $S_d$,
acting transitively on complete frames, so $S$ is connected. Reindexing
frames gives the left action
\[
                         L_a(\phi)=\phi\circ a^{-1}.
\]
It commutes with every isomorphism of the original covers. This gives
a named $G$-action on families of frame surfaces, rather than an action
defined only up to inner automorphism.

Fill the finite punctures and the ends above the outside point to obtain
a compact oriented surface $\overline S$. It is the regular $G$-cover
of the sphere with inertia $\langle\tau\rangle$ at each of the $N$
simple branch points and inertia $\langle\sigma\rangle$ at infinity.
Replacing $\sigma$ by its inverse leaves the latter subgroup unchanged.
This includes $\sigma=1$, when infinity is unramified. Filling ends
preserves connectedness.

\begin{lemma}\label{lem:cw}
In the rational representation ring of $G$,
\begin{equation}\label{eq:cw}
 [H^1(\overline S)]
 =2[\one]+(N-1)[\Q G]
       -N[\Q(G/\langle\tau\rangle)]
       -[\Q(G/\langle\sigma\rangle)].
\end{equation}
The trivial representation does not occur. For a nontrivial irreducible
representation $V_\rho$, put
\[
 n_\rho=\dim V_\rho,\quad
 a_\rho=\codim V_\rho^\tau,\quad
 c_\rho=\codim V_\rho^{\langle\sigma\rangle}.
\]
Its multiplicity in $H^1(\overline S)$ is
\begin{equation}\label{eq:multiplicity}
 m_\rho=Na_\rho-2n_\rho+c_\rho,
 \qquad n_\rho\le r a_\rho.
\end{equation}
In particular all nontrivial types occur if $N>2r$, and $m_\rho\ge4$
for all of them if $N\ge2r+4$.
\end{lemma}
\begin{proof}
The sphere minus its $N+1$ distinguished points has Euler characteristic
$1-N$. Its equivariant cellular Euler characteristic upstairs is
$(1-N)[\Q G]$. Filling a puncture with inertia $H$ adds the permutation
module $\Q(G/H)$. Both $H^0(\overline S)$ and $H^2(\overline S)$
are trivial, since the surface is connected and the deck action preserves
orientation. This proves \eqref{eq:cw}. The trivial multiplicity is
$2+(N-1)-N-1=0$. Frobenius reciprocity gives
\[
 (N-1)n_\rho-N(n_\rho-a_\rho)-(n_\rho-c_\rho)
                  =Na_\rho-2n_\rho+c_\rho.
\]
The $r$ adjacent transpositions generate $S_d$. Their fixed spaces in
a nontrivial irreducible have zero intersection, and each has codimension
$a_\rho$. Codimensions are subadditive under intersection, giving
$n_\rho\le ra_\rho$. Also $a_\rho\ge1$; otherwise all transpositions
would act trivially. Thus $m_\rho\ge(N-2r)a_\rho+c_\rho$, proving the
thresholds, without any exception for one-dimensional representations.
\end{proof}

Put $H=H^1(\overline S;\Q)$ and $K=H^1(S;\Q)$. The puncture sequence
embeds $H$ into $K$, and $H^2(S)=0$. The field $\Q$ is a splitting
field for $S_d$: its irreducible rational representations are the
absolutely irreducible Specht modules. Each admits a nondegenerate
invariant symmetric form, obtained by averaging a positive definite
rational form. Hence
\begin{equation}\label{eq:symplectic-decomp}
 H=\bigoplus_{\rho\ne1} V_\rho\otimes W_\rho,
 \qquad \langle-,-\rangle_H
       =\bigoplus_{\rho\ne1} B_\rho\otimes\omega_\rho,
\end{equation}
where $B_\rho$ is symmetric nondegenerate and $\omega_\rho$ is alternating
nondegenerate. Schur's lemma proves this decomposition of the compact
intersection pairing. In particular $\dim W_\rho=m_\rho$ is even.
Semisimplicity permits an equivariant splitting $K=H\oplus Z$. This is
a splitting of vector spaces for a fixed surface; no flat splitting in
a family will be used.

\subsection{The orbit-diagonal Leray sequence}
The space $O_m(S)$ is obtained from $S^m$ by removing the labelled
diagonals $y_i=g y_j$. A small loop around one of these diagonals gives
a meridian class. The Leray sequence organizes cohomology by the number
of such meridians and by the cohomology carried along their intersection
strata. In the degrees used below, only single diagonals and intersections
of two diagonals occur.

Consider the inclusion $j:O_m(S)\hookrightarrow S^m$, for $m=2,3$.
Its Leray sequence has
\[
 E_2^{p,q}=H^p(S^m;R^qj_*\Q)\Longrightarrow H^{p+q}(O_m(S);\Q).
\]
For $i<j$ and $a\in G$, write
\[
 \Gamma_{ij}(a)=\{y_i=a y_j\},\qquad E_{ij}(a)\in E_2^{0,1}.
\]
The latter symbol denotes the positive normal meridian. The action on
cohomology is inverse pullback, so
\begin{equation}\label{eq:edge-action}
 (a_1,\ldots,a_m)E_{ij}(g)=E_{ij}(a_i g a_j^{-1}).
\end{equation}

We recall the local construction to specify the differential we need.
Near a diagonal stratum, choose small disks at its distinct $G$-orbits
and identify their translates by the specified deck maps. The excluded
equations become the usual diagonals in an oriented plane. Local
meridian cohomology gives the Orlik--Solomon relations; orientations
trivialize the normal local systems. Along a fixed labelled stratum
the labels do not permute. The resulting coefficient systems are
constant. This is the orbit version of the construction in
\cite[Theorem 1]{Totaro} and \cite[proof of Theorem 3.1]{Casto}.

In the degrees needed here,
\begin{equation}\label{eq:leray-low}
 E_2^{p,0}=H^p(S^m),\qquad
 E_2^{p,1}=\bigoplus_{i<j,\,g\in G} H^p(\Gamma_{ij}(g)).
\end{equation}
For $m=3$, $R^2j_*\Q$ is a direct sum over labelled triple diagonals
of their constant rank-two local Orlik--Solomon spaces. Inconsistent
parallel-edge intersections are empty, because the action on $S$ is
free. The relations include
\begin{align}
 E_{ji}(g^{-1})&=E_{ij}(g),& E_{ij}(g)^2&=0,\notag\\
 E_{ij}(g)E_{ij}(h)&=0\quad(g\ne h),\label{eq:os-rel}\\
 E_{12}(g)E_{23}(k)
  &=E_{12}(g)E_{13}(gk)+E_{13}(gk)E_{23}(k).\notag
\end{align}
Total degrees anticommute. Coefficients restrict to their graph:
\begin{equation}\label{eq:coefficient-relation}
           (p_i^*\alpha-p_j^*g^*\alpha)E_{ij}(g)=0.
\end{equation}
Thus a forest is not a free module over all of $H^*(S^m)$.

Orient the normal meridian so that the single-divisor transgression is
the Thom class. Naturality and multiplicativity give
\begin{equation}\label{eq:actual-differential}
 d_2E_{ij}(g)=\delta_{ij}(g)=[\Gamma_{ij}(g)],
 \qquad |d_2|=(2,-1).
\end{equation}
For ordinary diagonals this is \cite[Theorem 2]{Totaro}, applied in
real dimension two. The differential is a derivation for total degree.
All constructions are $G^m$-equivariant. We use this local description
of the sequence and compute its differentials; no degeneration statement
for the noncompact surface is required.

\subsection{The graph-class kernel}
\begin{lemma}\label{lem:graph-kernel}
If $N>2r$, the map
\[
 D:\Q[G]\longrightarrow K\otimes K,\qquad[g]\longmapsto\delta(g),
\]
has kernel $\Q\nu$, where $\nu=\sum_{g\in G}[g]$.
\end{lemma}
\begin{proof}
The compact graph in $\overline S^2$ has two point-class terms and a
mixed term equal, up to one common orientation sign, to
$(g\otimes1)\Omega_H$, the transform of the compact symplectic
coevaluation. This follows by pairing its class with external products
and integrating their restrictions to the graph. Restriction to $S^2$
kills both point terms, since $H^2(S)=0$, and loses none of the mixed
term, since $H\hookrightarrow K$ is injective. Gysin classes commute
with this restriction. Therefore
\[
                         \delta(g)=(g\otimes1)\Omega_H
\]
with that fixed sign convention.

Identify $H\otimes H$ with $\End(H)$ by the compact pairing. The graph
map becomes the group-algebra action on $H$. Since
$\Q G=\bigoplus_\rho\End(V_\rho)$, \cref{lem:cw} says that its only
undetected block is the trivial block, spanned by $\nu$. The inclusion
in $K\otimes K$ does not change the kernel. The argument allows every
boundary product $\sigma$.
\end{proof}

\subsection{The two-coordinate theorem}
\begin{proposition}\label{prop:p2}
For $N>2r$, $P_2H^b(O_2(S);\Q)=0$ for $b\le1$.
\end{proposition}
\begin{proof}
The degree-zero term and $E_2^{1,0}=K\oplus K$ have zero $P_2$ part.
The other total-degree-one term is $E_2^{0,1}=\Q G$. Its differential
is $D$, whose kernel is the norm line. Averaging either coordinate
fixes this line, so $P_2$ kills it. The projected pieces vanish on $E_3$.
The exact idempotent functor and the finite abutment filtration give
the assertion.
\end{proof}

\subsection{The three-position differential in bidegree \texorpdfstring{$(1,1)$}{(1,1)}}
For three coordinates, counting the coordinates appearing in a class
leaves two possible contributions in total degree two. The first is a
meridian multiplied by a degree-one class on the remaining coordinate;
the second is a product of two meridians. We treat them in that order.
For the first contribution, the compact diagonal class supplies an
alternating tensor in two positions. The following calculation explains
why contributions from the three choices of positions cannot cancel.

On a graph in $S^3$, its degree-one coefficient is the sum of $K$ on
the merged coordinate and $K$ on the unused coordinate. The first
summand is killed by averaging the unused coordinate. Hence
\begin{equation}\label{eq:11-source}
 P_3E_2^{1,1}
 =P_3\bigl((\Q G)_{12}\otimes K_3
       \oplus(\Q G)_{13}\otimes K_2
       \oplus(\Q G)_{23}\otimes K_1\bigr).
\end{equation}
The target $E_2^{3,0}$ is $K_1\otimes K_2\otimes K_3$.

For a source element $E_{ij}(g)\alpha_k$, the differential is
$\delta_{ij}(g)\alpha_k$. Ordering the three odd cohomology factors by
coordinates gives signs $(+,-,+)$ for edges $(12),(13),(23)$.

\begin{lemma}\label{lem:three-position}
Let $(W,\omega)$ be symplectic of dimension $m\ge4$, and let
$\Omega$ be its alternating coevaluation. For any vector-space
extension $U=W\oplus Z$, the map
\begin{equation}\label{eq:three-map}
 U^{\oplus3}\longrightarrow U^{\otimes3},\qquad
 (a,b,c)\longmapsto\Omega_{12}a_3-\Omega_{13}b_2+\Omega_{23}c_1
\end{equation}
is injective.
\end{lemma}
\begin{proof}
Change the middle source coordinate by $b\mapsto-b$. On $W^{\oplus3}$,
contract the resulting plus-sign map in positions $12,13,23$.
With $\Omega=\sum_t(e_t\otimes f_t-f_t\otimes e_t)$ and
$\omega(e_t,f_t)=1$, the contractions are
\[
 \begin{pmatrix}m&1&-1\\1&m&1\\-1&1&m\end{pmatrix}
 \begin{pmatrix}a\\b\\c\end{pmatrix}.
\]
The determinant is $(m-2)(m+1)^2$, which is nonzero. If the unpaired
vector lies in $Z$, the three images lie in the separate direct summands
$W\otimes W\otimes Z$, $W\otimes Z\otimes W$ and
$Z\otimes W\otimes W$. Each map is injective and none can cancel a
term in $W^{\otimes3}$.
\end{proof}

We now identify this tensor map with \eqref{eq:actual-differential},
including all representation factors. Write
\[
 K=\bigoplus_\rho V_\rho\otimes U_\rho,
 \qquad U_\rho=W_\rho\oplus Z_\rho.
\]
The $\rho$ block of the graph map is
\[
 \End(V_\rho)\longrightarrow
 (V_\rho\otimes V_\rho)\otimes(W_\rho\otimes W_\rho),
 \quad T\longmapsto(T\otimes1)B_\rho^{-1}\otimes\Omega_\rho.
\]
The first factor is an isomorphism onto $V_\rho\otimes V_\rho$;
there is no unaccounted multiplicity or quotient. Only nontrivial
coordinate types survive $P_3$. An edge contributes to a type
$(\rho_1,\rho_2,\rho_3)$ precisely when its two endpoint types agree.
If all three types are not equal, at most one edge contributes, and
its map is injective. Different types belong to different $G^3$
isotypic summands.

If all three types equal $\rho\ne1$, factor out
$V_\rho^{\otimes3}$. The remaining differential is exactly
\eqref{eq:three-map} on $U_\rho^{\oplus3}$. When $N\ge2r+4$,
\cref{lem:cw} gives $\dim W_\rho\ge4$. Thus $d_2$ is injective on
$P_3E_2^{1,1}$, including every puncture-cohomology summand.

\subsection{The forest differential in bidegree \texorpdfstring{$(0,2)$}{(0,2)}}
A labelled triple diagonal is uniquely described by
$y_1=g y_2=h y_3$; the third edge has label $k=g^{-1}h$.
Its two independent local degree-two classes can be chosen as
\[
 A(g,h)=E_{12}(g)E_{13}(h),\qquad
 B(g,k)=E_{12}(g)E_{23}(k).
\]
Indeed the third product $C(h,k)=E_{13}(h)E_{23}(k)$ is related by
$A-B+C=0$, and the ordinary three-point local Orlik--Solomon space
has dimension two. The sheaf decomposition above shows that there
are no relations between different labelled triple diagonals. Their
degree-zero coefficient is $\Q$, since $S$ is connected. Consequently
\begin{equation}\label{eq:forest-basis}
                  E_2^{0,2}=\Q[G^2]_A\oplus\Q[G^2]_B.
\end{equation}
Both summands are $G^3$-stable. Their actions are respectively
\[
 (g,h)\mapsto(a_1ga_2^{-1},a_1ha_3^{-1}),\qquad
 (g,k)\mapsto(a_1ga_2^{-1},a_2ka_3^{-1}).
\]
Thus each is the permutation module on $G^3/G_{\rm diag}$.

The derivation formula gives
\begin{align*}
 d_2A(g,h)&=\delta_{12}(g)E_{13}(h)-E_{12}(g)\delta_{13}(h),\\
 d_2B(g,k)&=\delta_{12}(g)E_{23}(k)-E_{12}(g)\delta_{23}(k).
\end{align*}
Project the target $E_2^{2,1}$ to the component with residual edge
$13$ labelled $h$. Only the first term of $d_2A(g,h)$ contributes.
Parametrize that graph by
\[
            (y_1,y_2)\longmapsto(y_1,y_2,h^{-1}y_1).
\]
After coefficient restriction the map is literally $[g]\mapsto D(g)$
in $H^2(S\times S)=K\otimes K$. Its kernel, for fixed $h$, is the
norm line. A general kernel vector is therefore
$\sum_h c_h\sum_g A(g,h)$, invariant under the entire coordinate-$2$
group. Its $P_3$ part is zero.

Projection to residual edge $23$ labelled $k$ similarly sees only
$B(g,k)$, again with coefficient $D(g)$, using the parametrization
$(y_1,y_2)\mapsto(y_1,y_2,k^{-1}y_2)$. Its kernel is invariant under
coordinate $1$, and also has zero $P_3$ part. The two chosen target
components preclude cancellation between the two forest families.
Thus $d_2$ is injective on $P_3E_2^{0,2}$ whenever $N>2r$.

\begin{theorem}[Three-coordinate projector]\label{thm:p3}
If $N\ge2r+4$, then
\[
                       P_3H^b(O_3(S);\Q)=0\qquad(b\le2).
\]
For any connected frame surface, the same vanishing holds for $b\le1$.
\end{theorem}
\begin{proof}
Every term of total degree at most one uses at most two coordinates,
so $P_3$ kills it. In total degree two, $E_2^{2,0}$ also uses at most
two coordinates. The other terms are $(1,1)$ and $(0,2)$, whose
projected differentials were just proved injective. Their $P_3$ parts
vanish on $E_3$. Every later page is a subquotient in an
equivariant spectral sequence. Since the image of a finite-group
rational idempotent is an exact functor, it also kills the finite
filtration of the abutment. Duality gives the corresponding homology
statements.
\end{proof}

\section{Insertion operators and the local collision calculation}
\label{sec:insertions}

The projector theorem concerns configurations on a fixed frame surface.
To use it for stabilization, we let the cover vary and replace each
marked unbranched disk by a disk containing an equal pair of branch
points. Changing its frame changes the transposition inserted there.
Differences of frame choices therefore realize the operators
$x_e-x_f$.

We carry out this construction with its right module structure intact.
This will allow us to tensor the insertion maps with $C_z$ and obtain
information about $V_2$, where a statement about $H_2(M)$ alone would
not suffice. The last part of the section studies what happens when
two marked disks approach one another. That collision gives the
stronger bound on $V_1$ used in the genus estimate.

\subsection{A family with a global frame action}
Let $E_m\to\CHur^{T_d}$ parametrize $m$ ordered unbranched interior
points, a complete sheet frame at each, and pairwise disjoint small
round disks about them. Remember the disk radii as part of the object.
At fixed points their allowed set is nonempty and convex, so forgetting
radii is an equivariant homotopy equivalence. Over a cover $f$ the fibre
is $O_m(S_f)$. Frame reindexing commutes with all transition maps and
defines an action of $G^m$ over the identity on the base.

The Serre sequence retains the possibly nonconstant coefficient systems
$\mathcal H_b(O_m(S_f))$. If $P_m$ kills those systems in degrees
$b\le i$, it kills every term contributing to total homological degree
at most $i$. Exactness of the idempotent and the finite filtration then
give $P_mH_{\le i}(E_m)=0$. In particular the preceding section proves
\begin{align}
 P_2H_1(E_2)_N&=0 &&(N>2r),\label{eq:family-p2}\\
 P_3H_2(E_3)_N&=0 &&(N\ge2r+4),\label{eq:family-p3}\\
 P_3H_1(E_3)_N&=0 &&\text{for every }N.\label{eq:family-p3-one}
\end{align}
The pointwise splitting $K=H\oplus Z$ in the preceding proof was used
only to prove fibrewise zero statements, not to trivialize these local
systems.

\subsection{Surgery and the right-module identity}
In each marked disk, the frame trivializes the old inverse image as
$d$ disks. Replace this trivial cover by the standard two-branch cover
with local tuple $((12),(12))$. Its boundary monodromy is one, so it
glues to the unchanged exterior cover. The local active component is
an annulus and the other sheets are disks. Removing unbranched disks
does not disconnect the old source; attaching these pieces preserves
connectedness. The exterior framing and boundary product are unchanged.
This defines
\[
                         F_m:E_m\longrightarrow\CHur^{T_d},
\]
raising word length by $2m$.

For precision, continuity can be checked in simply connected
configuration charts. The branch points and remembered disks move by
ambient isotopy, the frames are discrete lifts, and the glued monodromy
is constant under this transport. Hence the construction lifts
continuously to the finite framed Hurwitz covering. On overlaps, the
cover isomorphism preserves the local frames and the gluing, so the
lifts agree.

Use Moore rectangles with the lower boundary path transporting the
numbering to the right gluing interval. Appending an arbitrary cover on
the right leaves all old marked disks and their radii unchanged. Thus
$E_m$ is a right module, the $G^m$-action is right linear, and $F_m$
is right linear. Prefix an empty fixed rectangle containing the $m$
disks, with frames transported from its left boundary. This defines
a right-module map $s_m:M\to E_m$; forgetting the empty prefix is
homotopic to the identity. The disks are ordered once and for all.

The placement of the disks explains the side of the module structure.
For insertions obtained from $s_m$, surgery takes place in a fixed
region to the left of the old cover,
whereas the right action appends a new cover beyond the old one. These
operations take place in disjoint regions, so their compatibility is
already present in the family of covers before passing to homology.

For $g=(g_1,\ldots,g_m)$, the composite $F_mL_gs_m$ is left insertion
of the ordered packets $(\tau^{g_j},\tau^{g_j})$. Any adjustment to
standard packet positions takes place in the fixed prefix and commutes
with right multiplication. For given $e_j,f_j\in T_d$, choose
$g_j,h_j$ with $\tau^{g_j}=e_j$ and $\tau^{h_j}=f_j$. Put
\[
 Q_m=\prod_j(L_{g_j}-L_{h_j}),\qquad Q_m=Q_mP_m.
\]
In the derived category of right $A$-modules, this gives
\begin{equation}\label{eq:operator-identity}
 F_mQ_ms_m=L_{\prod_j(x_{e_j}-x_{f_j})}:M\sh{2m}\longrightarrow M.
\end{equation}
The order of the packets is fixed during this identity. Centrality in
$R$ subsequently extends it to the ideal action.

\begin{proposition}\label{prop:high}
The square operators satisfy
\begin{align}
 \Delta^3H_1(M)&=0,\label{eq:abs-global-one}\\
 \Delta^2H_1(M)_N&=0 &&(N>2r),\label{eq:abs-high-one}\\
 \Delta^3H_2(M)_N&=0 &&(N\ge2r+4).\label{eq:abs-high-two}
\end{align}
For the derived cofiber $X=M\dt_A C_z$,
\begin{equation}\label{eq:relative-high}
                 \Delta^3V_{2,L}=0\qquad(L\ge2r+4).
\end{equation}
\end{proposition}
\begin{proof}
The absolute statements follow by applying
\eqref{eq:family-p2}--\eqref{eq:family-p3-one} to
\eqref{eq:operator-identity}. To prove the relative statement, tensor
the right-module identity with the fixed bimodule $C_z$. There is a
$G^3$-equivariant triangle
\[
 C_*(E_3)\sh2\longrightarrow C_*(E_3)
                    \longrightarrow C_*(E_3)\dt_A C_z.
\]
In homological degree two and word $L$ its outer controlling terms are
$H_2(E_3)_L$ and $H_1(E_3)_{L-2}$. Their $P_3$ parts vanish by
\eqref{eq:family-p3} and \eqref{eq:family-p3-one}. Applying the exact
idempotent to the long exact sequence proves the relative vanishing.
Thus the relative assertion uses the degree-two projector theorem at
word $L$, together with the all-word degree-one assertion at word
$L-2$. The tensor product throughout is taken with $C_z$ on the right.
\end{proof}

\begin{corollary}\label{cor:j-high}
The following are source-word annihilation statements:
\begin{align}
 J^3V_{2,L}&=0 &&(L\ge2r+4),\label{eq:jv2}\\
 J^2\Kh_{1,L}&=0 &&(L\ge2r+3),\label{eq:jk1}\\
 J^3\Kh_{2,L}&=0 &&(L\ge2r+6).\label{eq:jk2}
\end{align}
Also $J^3V_1=0$ globally.
\end{corollary}
\begin{proof}
On $\Kh_i$ and $C_i$, the action of $z$ is zero, so $J=(z,\Delta)$
has the same powers as $\Delta$. The kernel $\Kh_{i,L}$ is a subspace
of $H_i(M)_{L-2}$; this explains the two-word shifts in
\eqref{eq:jk1} and \eqref{eq:jk2}. The global $V_1=C_1$ statement
follows from \eqref{eq:abs-global-one}.

For $L\ge2r+3$, $\Delta^2V_{2,L}$ maps to zero in
$\Kh_{1,L+4}$, hence lies in $C_{2,L+4}$ and is killed by $z$.
Thus $z\Delta^2V_{2,L}=0$ there. Expand each cubic in
$J=(z,\Delta)$. Terms with no $z$ vanish by \eqref{eq:relative-high},
terms with one $z$ by the preceding observation, and terms with at least
two $z$ by \eqref{eq:z-square}. This proves \eqref{eq:jv2}.
\end{proof}

\subsection{Generators of connected first homology}
\begin{lemma}\label{lem:h1-generation}
As a left $R$-module, $H_1(M)$ is generated in word degrees at most
$r+1$.
\end{lemma}
\begin{proof}
A connected Hurwitz component is the classifying space of the stabilizer
of a numbered transposition tuple. The liftable-braid theorem of
Mulazzani--Piergallini \cite[Theorem B and Lemmas 2.1, 2.3]{MP} says
that this group is generated by liftable powers of half twists along
intervals. The boundary-fixing liftable group is the numbered stabilizer:
fixing the boundary fibre fixes each boundary circle by unique lifting.
Their theorem allows arbitrary boundary monodromy.

Transport the interval of such a generator to the final adjacent pair.
The resulting local loop is $\sigma_{N-1}^k$, where the minimal positive
$k$ is $1$, $2$ or $3$ for equal, disjoint or overlapping transpositions.
The transport braid is a path in the action groupoid; it need not itself
be a stabilizer. Change of basepoint transports its loop class in $H_1$.

Let $s$ be the support rank of the final pair: $s=1$ for equal labels
and $s=2$ otherwise. If the current tail is not connected, a letter
of the remaining prefix crosses two of its support blocks. Move that
letter to the end of the prefix by inverse Hurwitz moves
$(a,c)\mapsto(aca,a)$, and add it to the tail. This merges exactly two
blocks. After $r-s$ steps the tail is connected and has length
$r-s+2\le r+1$.

This extraction uses only the first $N-2$ positions. Its braid lies in
$B_{N-2}$ and commutes with $\sigma_{N-1}^k$. The resulting naturality
square in the action groupoid is a cylinder between the local loops.
It identifies the two loop classes in the homology of the connected
Hurwitz space.
The final class is the product of a prefix coefficient in $R$ and an
$H_1$ class on the connected tail. Abelianizing the interval generators
proves the lemma.
\end{proof}

\subsection{A collision produces a local four-branch loop}
The next step is independent of the three-coordinate calculation.

\begin{lemma}\label{lem:local-four}
For a four-branch tuple $(a,a,b,b)$, the first homology of its local
Hurwitz component is generated by the image of
\[
                            R_2\otimes H_1(A)_2.
\]
All three possibilities for $a,b$ are included.
\end{lemma}
\begin{proof}
If $a=b$, the active two-sheet cover is connected. If $a,b$ overlap,
the active three-sheet cover is connected. In both cases apply the
interval-power theorem of \cite{MP} to the active cover; unused sheets
are fixed unramified disks. Move a generating interval to the final
pair. The loop is then a two-letter coefficient followed by a
two-letter loop. Transport in the action groupoid identifies the
homology classes, as in \cref{lem:h1-generation}.

If $a,b$ are disjoint, their conjugations commute and the component is
the configuration space $\Conf_{2,2}(D^\circ)$ of two points of each
colour. It is the complement of the cross-colour collision divisor in
$\UConf_2(D^\circ)\times\UConf_2(D^\circ)$. Its $H_1$ is generated
by the two within-colour half twists and the cross-colour meridians.
Indeed 1-cycles in the ambient product can avoid the divisor, and a
2-chain bounding a nullhomologous cycle can be made transverse to its
smooth part and disjoint from its multiple-collision strata, which
have real codimension at least four. Excision expresses its boundary
as a sum of meridians. A cross-colour meridian is a squared half twist
of the corresponding disjoint pair. Transporting each supporting
interval to the final pair proves the same factorization. Singleton
sheets contribute no additional motion.
\end{proof}

\begin{proposition}[Collision divisibility]\label{prop:collision}
For source word degree $N\ge2r-1$,
\begin{equation}\label{eq:collision}
 \Delta^2H_1(M)_N\subseteq zH_1(M)_{N+2}
                      \subseteq H_1(M)_{N+4}.
\end{equation}
\end{proposition}
\begin{proof}
There are two steps. First, a collision turns the projected motion of
the markers into a loop involving only four new branch points. Second,
we factor that local loop so that two of its letters join the old
connected cover as an $H_0$ coefficient. It is this connected
coefficient that we will divide by $z$.

Fix a base component $B_0\subset\CHur_N$ and its two-marker family,
with connected fibre $F=O_2(S)$. The five-term homology sequence contains
\[
 H_2(B_0;\Q)\longrightarrow H_0(B_0;\mathcal H_1(F))
 \longrightarrow H_1(E_2)\longrightarrow H_1(B_0;\Q)\longrightarrow0.
\]
The two frame actions commute with monodromy and act trivially on the
base and $\mathcal H_0(F)$. Applying $P_2$ shows that projected fibre
classes generate $P_2H_1(E_2)$. This uses the local
system and its coinvariants into account.

The graphs $\Gamma_g\subset S^2$ are pairwise disjoint. Excision for
their oriented normal disks gives
\[
 H_2(S^2)\longrightarrow\bigoplus_{g\in G}H_0(\Gamma_g)
 \longrightarrow H_1(O_2(S))\longrightarrow H_1(S^2)\longrightarrow0.
\]
Since $P_2H_1(S^2)=0$, projected diagonal meridians generate the reduced
fibre group. This is ordinary relative homology; noncompactness of
the graphs adds no degree-one generator.

A small graph meridian moves one marker around the other in an old
unbranched disk. The old cover trivializes there and both frames return
unchanged. Under surgery only the four new branch points move, with
local tuple $(a,a,b,b)$; the exterior cover remains fixed. One model
is the square of the packet exchange
$\beta=\sigma_2\sigma_1\sigma_3\sigma_2$. A tether from this disk to
the external boundary separates it from a complementary disk containing
all old branch values. Removing that unbranched disk and its tether
cuts only trivial strips, so the old cover remains connected. Isotopy
to the concatenation model is a groupoid transport of the loop.
Passing the packet through the old cover can simultaneously conjugate
its labels by the old boundary product, but preserves the three support
types of \cref{lem:local-four}.

It follows from that lemma and \eqref{eq:operator-identity} that
\[
                  \Delta^2H_1(M)_N\subseteq
                  H_0(M)_N\,R_2\,H_1(A)_2.
\]
The coefficient lies in \emph{connected} $H_0(M)_{N+2}$, because the
old tuple already generates $S_d$. In each product sector
$\sigma'$, \eqref{eq:hzero} starts in degree $2r-\ell(\sigma')\le2r$.
Since $N+2>2r$, the coefficient contains a positive power of $t$ and
is divisible by $z=qt$, leaving a connected coefficient in degree $N$
with the same product. Associativity gives
$(zb')\lambda=z(b'\lambda)$, with $b'\lambda\in H_1(M)_{N+2}$.
The word degrees can be read directly from this factorization. The
old coefficient has degree $N$, the extra two letters raise it to
$N+2$, and division by $z$ returns it to degree $N$. Multiplication
by the remaining two-letter loop then gives $b'\lambda$ in degree
$N+2$, whose $z$-multiple is the original class in degree $N+4$.
This proves \eqref{eq:collision}. Its image becomes zero in $V_1=C_1$,
which is the consequence needed next.
\end{proof}

\begin{corollary}\label{cor:h1-new}
The first cofiber homology satisfies $\Endw V_1\le2r+2$.
\end{corollary}
\begin{proof}
We know $V_1=C_1$, $zV_1=0$, $J^3V_1=0$, and generation in words
at most $r+1$. The collision statement gives
$\Delta^2V_{1,N}=0$ for $N\ge2r-1$.
Consider a homogeneous term in degree $L\ge2r+3$. Its coefficient
on a chosen generator has degree $q'\ge r+2$, and by \cref{lem:ring}
is a sum of terms $x_{e_1}\cdots x_{e_k}u_\sigma$, where
$q'=2k+\ell(\sigma)$. Terms with $k\ge3$ vanish, and $k=0$ is
impossible. For $k=1$, necessarily $q'=r+2$ and $\ell(\sigma)=r$.
All squares agree in this top product sector, so
$x_eu_\sigma=q^{-1}zu_\sigma$ and the term vanishes. For $k=2$ the
source degree is $L-4\ge2r-1$; since $zV_1=0$, this two-square action
equals a $\Delta^2$ action and vanishes by collision divisibility.
Thus $V_{1,L}=0$ for every $L\ge2r+3$.
\end{proof}

\section{Kernel bounds and the uniform genus range}
\label{sec:bounds}

We now assemble the preceding estimates. The quotient by $J$ is
controlled by the geodesic coefficient calculation, while powers of $J$
vanish above the source degrees found in the insertion argument. A
short lifting lemma combines these two kinds of information. Starting
with $V_0$ and the collision bound for $V_1$, we then bound $V_2$ and
the first two stabilization kernels. The final step translates the
kernel indices into genus and passes from cohomology to line bundles.

\subsection{Lifting a bound with localized nilpotence}
\begin{lemma}\label{lem:local-lift}
Let $J$ be generated by central degree-two elements of a nonnegatively
word-graded algebra. Suppose
\[
 \Endw(V/JV)\le\beta,\qquad
 J^aV_L=0\quad\text{for every source word }L\ge T.
\]
Then
\begin{equation}\label{eq:local-lift}
             \Endw V\le\max\{\beta+2a-2,\ T+2a-1\}.
\end{equation}
\end{lemma}
\begin{proof}
For $n>\beta$, $V_n=\sum_e x_eV_{n-2}$. If
$n>\beta+2a-2$, repeat this equality $a$ times; its final source degree
is $n-2a$. If also $n>T+2a-1$, integrality implies $n-2a\ge T$,
so the resulting $a$-fold products vanish. Centrality absorbs all
algebra coefficients into the module. This proves the strict endpoints.
\end{proof}

The second endpoint is necessary in general. The module
$\Q[x]\sh{T-1}/(x^{a+1})$, with $|x|=2$, is killed by $x^a$ on
all source degrees at least $T$, but its top nonzero degree is
$T+2a-1$. A high-word annihilation hypothesis therefore cannot be
silently replaced by global nilpotence.

\subsection{The relative Postnikov calculation}
Set $F(-)=B\dt_A(-)$ and $Y=B\dt_A M$, retaining its right $A$-action.
Associativity, with the original factor order, gives
\[
 F(X)\simeq (B\dt_A M)\dt_A C_z.
\]
The cofiber triangle $Y\sh2\to Y\to F(X)$ and \eqref{eq:bcm} imply
\begin{equation}\label{eq:fx-bound}
                     \Endw H_iF(X)\le2r+i+1\quad(i\ge1).
\end{equation}
Indeed, the two terms controlling $H_iF(X)$ are $H_i(Y)$ and
$H_{i-1}(Y)\sh2$; the latter accounts for the extra one in the bound.

The next lemma measures the contribution of the lower homology of $X$
to $V_i/JV_i$. A Postnikov tower separates $X$ into its homology
modules. Applying the coefficient functor to that tower need not be
exact on each module, and the terms one homological degree higher are
what enter the estimate.

\begin{lemma}\label{lem:bpost}
If $h_j\ge\Endw V_j$ for $j<i$, then
\begin{equation}\label{eq:bpost}
 \beta_i:=\Endw(V_i/JV_i)
 \le\max\left\{2r+i+1,\ \max_{0\le j<i}
                           \bigl(h_j+2(i-j+1)\bigr)\right\}.
\end{equation}
\end{lemma}
\begin{proof}
The connective functor $F$ preserves lower connectivity, so
$H_iF(X)=H_iF(\tau_{\le i}X)$. Apply it to the Postnikov triangle
\[
 V_i[i]\longrightarrow\tau_{\le i}X\longrightarrow
 \tau_{\le i-1}X\longrightarrow V_i[i+1].
\]
Its exact sequence includes
\[
 H_{i+1}F(\tau_{\le i-1}X)\longrightarrow
 H_0F(V_i)\longrightarrow H_iF(X).
\]
Since $A$ is connective and $B,V_i$ are discrete,
$H_0F(V_i)=B\otimes_RV_i=V_i/JV_i$.
The finite lower Postnikov tower has pieces $V_j[j]$. Each $V_j$ has
nonnegative word grading, so the lower-support hypothesis of
\cref{cor:b-cost} is satisfied. Its incoming term has word support at most
$\max_{j<i}(h_j+2(i-j+1))$. Combining this with \eqref{eq:fx-bound}
proves the result. The relevant Tor degree is $i-j+1$.
\end{proof}

\begin{proposition}\label{prop:kernel-bounds}
For every $d\ge6$, with $r=d-1$,
\begin{align*}
 \Endw V_0&\le2r,& \Endw V_1&\le2r+2,&
 \Endw V_2&\le2r+10,\\
 \Endw\Kh_1&\le2r+8,& \Endw\Kh_2&\le2r+18.
\end{align*}
\end{proposition}
\begin{proof}
The first two inequalities are \eqref{eq:low-cone} and
\cref{cor:h1-new}. The next Postnikov step gives
\[
 \beta_2\le\max\{2r+3,2r+6,2r+6\}=2r+6.
\]
Apply \cref{lem:local-lift} to \eqref{eq:jv2}, with $a=3$ and
$T=2r+4$. It gives
\[
 h_2\le\max\{2r+6+4,2r+4+5\}=2r+10.
\]
Now
\[
 \beta_3\le\max\{2r+4,2r+8,2r+8,2r+14\}=2r+14.
\]
The surjection $V_3\twoheadrightarrow\Kh_2$ induces a surjection
$V_3/JV_3\twoheadrightarrow\Kh_2/J\Kh_2$. Thus its quotient modulo
$J$ has end at most $\beta_3$. Applying \cref{lem:local-lift} to
\eqref{eq:jk2}, with $a=3,T=2r+6$, gives
\[
 \Endw\Kh_2\le\max\{2r+14+4,2r+6+5\}=2r+18.
\]
Similarly $V_2\twoheadrightarrow\Kh_1$ passes the bound $\beta_2$
to its quotient modulo $J$; use $a=2,T=2r+3$ to obtain
\[
 \Endw\Kh_1\le\max\{2r+6+2,2r+3+3\}=2r+8.
\]
The surjection from $V_3$ is used only after taking the quotient by $J$.
The localized annihilation statement for $\Kh_2$ then finishes the
estimate, so a bound on all of $V_3$ is unnecessary.
\end{proof}

\subsection{From injectivity to the fixed-target space}
Now restrict to the connected product-one component. A cover of the
sphere of genus $g$ has
\[
                         N=2g+2d-2=2g+2r
\]
simple branch points. The kernel of stabilization with source word
$N$ is $\Kh_{i,N+2}$. By \cref{prop:kernel-bounds}, the maps in
degrees one and two are injective whenever
\begin{equation}\label{eq:genus-endpoint}
                         N+2>2r+18.
\end{equation}
For example, when $d=6$ we have $r=5$. At $g=8$ the source word is
$N=26$, so the kernel is indexed by $N+2=28$, exactly the endpoint
$2r+18$. At $g=9$ its index is $30$, strictly beyond that endpoint.
The general calculation is the same, since $N+2-(2r+18)=2g-16$.

All subsequent two-word stabilizations stay in this range. Iterating
to the eventual componentwise comparison of
\cite[Theorems 1.4.1 and 1.4.2]{LL} gives
$H_2=0$ and $\dim H_1\le1$. The branch map is a finite covering
of the unordered affine configuration space. Transfer makes the pullback
of its discriminant class nonzero in $H^1$, so $H_1=\Q$.

Removing the sheet numbering takes finite-group invariants. The
discriminant class is invariant, so the unnumbered affine Hurwitz space
still has $H^1=\Q$ and $H^2=0$.

Let $\mathcal H$ be the space with fixed target $\PP^1$, and let
$\mathcal H^{\rm aff}$ be the open where infinity is not a branch value.
Its complement $D_\infty$ is a smooth divisor: the branch map is
\'etale over the simple configuration locus. This divisor is connected.
Indeed, after removing the transposition at infinity, the remaining
tuple still generates $S_d$, since the removed transposition is the
inverse product of the others. The orbit theorem gives one component
for each fixed transposition product, and sheet conjugation permutes
these transitively.

The affine discriminant class has nonzero residue at this divisor.
For example, writing the ordinary discriminant as
$\prod_{i<j}(b_i-b_j)^2$, a branch coordinate tending to infinity
with local coordinate $t=1/b_1$ contributes order $-(2N-2)$.
The low-degree Gysin sequence is therefore
\[
 0\to H^1(\mathcal H)\to\Q
 \xrightarrow{\,-(2N-2)\,}\Q
 \to H^2(\mathcal H)\to0.
\]
It follows that $H^1(\mathcal H)=H^2(\mathcal H)=0$.
Taking the target quotient is the homotopy quotient by $\PGL_2$.
Since $H^j(B\PGL_2;\Q)=0$ for $1\le j\le3$, its spectral sequence
gives
\begin{equation}\label{eq:hurwitz-low-cohom}
                  H^1(\U_{d,g};\Q)=H^2(\U_{d,g};\Q)=0.
\end{equation}

\subsection{Why this implies the full Picard statement}
We record the algebraic interface explicitly, since topological
$H^2$-vanishing by itself is not a general statement about holomorphic
line bundles on arbitrary nonproper spaces.
For an algebraic line bundle, a smooth compactification gives a useful
description of the kernel of the Chern class: after subtracting boundary
divisors, the obstruction lies in the Picard variety of the
compactification. The $H^1$ hypothesis below eliminates its invariant
part. This is why both low-degree cohomology groups were computed.

\begin{lemma}\label{lem:pic-c1}
Let $\mathcal X=[X/\Gamma]$, where $X$ is a smooth quasiprojective
complex variety and $\Gamma$ is finite. If
$H^1(\mathcal X;\Q)=0$, then
\[
             c_1:\Pic(\mathcal X)_\Q\longrightarrow
                                     H^2(\mathcal X;\Q)
\]
is injective. Restriction to a smooth open substack is surjective on
rational Picard groups.
\end{lemma}
\begin{proof}
Choose a $\Gamma$-equivariant smooth projective compactification
$\overline X$ with normal-crossing boundary. Purity makes
$H^1(\overline X)\to H^1(X)$ injective, so
$H^1(\overline X)^\Gamma=0$. Hodge theory implies that the fixed
identity component of $\Pic^0(\overline X)$ has dimension zero.
Therefore $(\Pic^0(\overline X)\otimes\Q)^\Gamma=0$: after clearing
denominators and averaging, an invariant point lies in the finite fixed
group and vanishes rationally.

Extend a line bundle from $X$ by closing a rational divisor in the
regular variety $\overline X$. Average its class over $\Gamma$.
If its Chern class is zero on $X$, the averaged extension has Chern
class equal to a rational combination of boundary divisors, by purity
in degree two. Subtract an averaged such combination. The resulting
invariant class lies in $\Pic^0(\overline X)_\Q$ and is zero. Thus
the original line is zero rationally on $X$ and, by finite norm, on
$\mathcal X$. For open restriction, the same divisor-extension and
averaging argument applies. The obstruction to a finite-group
linearization is torsion and disappears over $\Q$.
\end{proof}

For $g\ge9$, a full source-curve level structure of level at least
three gives such a finite quotient presentation of $\U_{d,g}$; the
level pencil scheme is quasiprojective by \cref{lem:full-pencil}.
That presentation is proved algebraically, independently of the
cohomology calculation in this section. Apply \cref{lem:pic-c1} to
\eqref{eq:hurwitz-low-cohom}. This proves \cref{thm:uniform}, since
\eqref{eq:genus-endpoint} is $2g+2>18$, whose first integer solution is
$g=9$. At $g=8$ there is equality and this argument gives no conclusion.
Finite inertia also identifies rational stack and coarse Picard groups:
a common tensor power kills stabilizer characters and descends, while
pullback is rationally injective.

\section{Universal pencils and the cases with \texorpdfstring{$\rho\ge3$}{rho at least 3}}
\label{sec:pencils}

We now give the geometric argument. It is independent of the
topological part. Throughout this section $g\ge6$, $d\le g-1$, and
$\rho=2d-g-2\ge3$ unless a different hypothesis is specified.

A basepoint-free pencil determines a map to a rational curve, with no
choice of coordinates on the target. Allowing base points gives a larger
space that is proper over the moduli of source curves. We will compute
the Picard group of this larger space and then restrict to the simple
Hurwitz open. The advantage is that properness makes fibrewise
cohomology comparisons available, while the complement of the simple
open supplies relations among the resulting line bundles.

\subsection{The full pencil stack}
Before rigidification, a pencil is a smooth genus-$g$ curve
$\pi:\mathcal C\to S$, a degree-$d$ line bundle $L$, a rank-two
bundle $V$, and an evaluation map $\pi^*V\to L$ inducing an injection
on sections on every geometric fibre. Evaluation need not be surjective,
and the pencil need not be complete. Simultaneous scalar multiplication
on $L$ and $V$ is a central $\Gm$. Rigidify it to obtain
$\Pencils$, and write
\[
             \Pencils\xrightarrow{q}\Jac\xrightarrow{a}\Mg,
             \qquad p=aq,
\]
where $\Jac$ is the similarly rigidified universal degree-$d$ Jacobian.
All remaining stabilizers embed in the finite automorphism group of
the source curve. The fibre of $p$ over a specified curve is the
scheme $G^1_d(C)$, not its quotient by source automorphisms.

\begin{lemma}\label{lem:full-pencil}
For $g\ge2$ and $d>0$, the stack $\Pencils$ is smooth of dimension
$2d+2g-5$ whenever nonempty, and $p$ is representable and proper.
After a finite source-level cover it is a smooth quasiprojective scheme.
\end{lemma}
\begin{proof}
\'Etale locally on $\Jac$, choose a Poincar\'e line bundle and a
sufficiently positive relative divisor $D$. The full pencil scheme is
the zero scheme of
\[
 S_2\longrightarrow F\quad\text{on }\Gr(2,E),\qquad
 E=\pi_*L(D),\quad F=\pi_*(L(D)|_D).
\]
Here $S_2$ is the tautological rank-two subbundle on the Grassmannian.
The vanishing of its map to $F$ says exactly that its two sections of
$L(D)$ vanish on $D$, and hence are sections of $L$.
Tensoring a Poincar\'e bundle by a base line tensors $E$ and $S_2$
simultaneously, so the zero scheme descends. It is proper over the
relative Jacobian, hence proper over $\Mg$.

For smoothness, choose an ordered independent pair $s_0,s_1$.
Let $\mathcal D(L)$ denote the sheaf of first-order differential
operators on $L$ with scalar symbol. It describes infinitesimal
automorphisms of the pair $(C,L)$, and its first cohomology describes
deformations of that pair. Adding the two sections gives the deformation
complex
\[
 [\mathcal D(L)\longrightarrow L\oplus L],\qquad
 P\longmapsto(Ps_0,Ps_1),
 \qquad0\to\mathcal O_C\to\mathcal D(L)\to T_C\to0,
\]
in cohomological degrees zero and one. In a trivialization its matrix is
\[
                 \begin{pmatrix}s_0&\partial s_0\\s_1&\partial s_1\end{pmatrix}.
\]
The determinant is a nonzero Wronskian: $s_1/s_0$ is nonconstant and
the characteristic is zero. Thus the map is injective with torsion
cokernel $Q$ of length $2d+2g-2$. The complex is $Q[-1]$; its
obstruction space is $H^1(Q)=0$. To see that this also applies with
base points, one can use \v Cech gluing: a deformation of the curve
and line has transition cochain in $\mathcal D(L)$, a deformation of
the two sections has cochain in $L\oplus L$, and their compatibility
is exactly the total-complex equation. No basepoint-free map to a
target is used. The ordered-pair problem therefore has dimension
$2d+2g-2$. Changing basis and removing the scalar stabilizer subtracts
three, giving the asserted dimension.

The descended determinant line constructed below is relatively ample,
up to sign and positive power, on the Grassmannian fibres. A theta
polarization on the Jacobian makes $p$ projective after a source-level
cover. The base of that cover is quasiprojective, and its level
structure kills the remaining automorphisms. This gives the final
assertion.
\end{proof}

The open of basepoint-free pencils with simple ramification and distinct
branch values is precisely $\U_{d,g}$. The target is unparameterized;
in a family its descent may be a conic rather than a globally trivial
$\PP^1$-bundle.

\subsection{Tautological line bundles and their descent}
On the stack before scalar rigidification, the universal line $L$ and
the rank-two bundle $V$ are available. They need not descend separately:
a scalar acts on each with weight one. The combinations below are
chosen so that these weights cancel. Thus their Chern classes will
represent line bundles on the rigidified stack where we need them.

On the scalar gerbe write
\[
 l=c_1(L),\quad k=c_1(\omega_\pi),\quad\gamma=c_1(V),\quad
 \alpha=\pi_*(l^2),\quad\beta=\pi_*(lk),\quad
 \lambda=c_1(\pi_*\omega_\pi).
\]
The Deligne pairings define weight-zero lines
\[
 \mathcal A=\langle L,L\rangle(\det V)^{-d},\qquad
 \mathcal B=\langle L,\omega_\pi\rangle(\det V)^{-(g-1)}.
\]
They descend to $\Pencils$, with classes
\begin{equation}\label{eq:a0b0}
                   a_0=\alpha-d\gamma,\qquad
                   b_0=\beta-(g-1)\gamma.
\end{equation}
Indeed the pairing weights are $2d$ and $2g-2$, and the weight of
$\det V$ is two.

Set $u=d+g-1$. The bundle $W=\pi_*(\omega_\pi L)$ has rank $u$,
since $H^1(\omega_C L)=H^0(L^{-1})^*=0$. The line
$\mathcal T=(\det V)^u(\det W)^{-2}$ also descends, and GRR gives
\begin{equation}\label{eq:tau-eta}
 \tau=c_1(\mathcal T)=-2\lambda-a_0-b_0,
 \qquad \eta_0=d b_0-(g-1)a_0.
\end{equation}
The second class is pulled back from the rigidified Jacobian, represented
there by $\langle L,\omega\rangle^d\langle L,L\rangle^{-(g-1)}$.
The determinant of the change from $(a_0,b_0)$ to
$(\eta_0,\tau+2\lambda)$ is $u\ne0$.

\subsection{Lefschetz comparison on every fixed curve}
Fixing the source curve gives a pencil scheme inside a Grassmannian
bundle over its Jacobian. Fulton--Lazarsfeld's incidence theorem compares
this zero locus with the ambient bundle by controlling the cohomology
of its complement. The comparison range depends on $\rho$, and the
positivity used in it holds for every smooth source curve. In particular,
we can later apply proper base change even over curves whose pencil
scheme is singular or has larger than expected dimension.

\begin{lemma}\label{lem:fl}
Fix any smooth genus-$g$ curve and assume $d\le g-1$. There is a
Grassmannian presentation $Y_C=G^1_d(C)\subset X_C$ for which
\begin{equation}\label{eq:fl-relative}
                         H_i(X_C,Y_C;\Q)=0\qquad(i\le\rho).
\end{equation}
For $\rho\ge3$ the inclusion induces isomorphisms through degree two;
for $\rho=2$ it induces an $H^1$ isomorphism and an $H^2$ injection.
\end{lemma}
\begin{proof}
Choose $m\ge2g$ sufficiently large, a point $P_0$, and distinct points
$P_1,\ldots,P_{m-d}$ with divisor $D$. On $\Pic^0(C)$ let $\mathcal P$
be a normalized Poincar\'e bundle and put
\[
 E=\operatorname{pr}_*(\mathcal P(mP_0)),\qquad
 F=\operatorname{pr}_*(\mathcal P(mP_0)|_D).
\]
Their ranks are $e=m+1-g$ and $f=m-d$, so $e\le f$.
Translation by $\mathcal O_C(mP_0-D)$ identifies the kernel at a
point with $H^0(C,L)$ of degree $d$. Let $X_C=\Gr(2,E)$ and
$Y_C=Z(S_2\to F)$. The positivity hypothesis for this
presentation is \cite[Lemma 2.2]{FL}. With $k_0=e-2$,
\cite[Proposition 1.2 and Remark 1.9]{FL}, equivalently the incidence
formulation in \cite[\S1.1]{Debarre}, gives the complement vanishing
\[
 H^j(X_C\setminus Y_C;\Q)=0\qquad(j\ge g+2(f+k_0)).
\]
Lefschetz duality and $\dim X_C=g+2(e-2)$ give the relative bound
\[
 2\dim X_C-g-2(f+k_0)=g-2(f-k_0)=2d-g-2=\rho.
\]
The endpoint is inclusive. The long exact sequence proves the asserted
isomorphisms and injection. No smoothness or expected dimension of
$Y_C$ is required.
\end{proof}

Grassmannian-bundle cohomology now identifies, for $\rho\ge3$,
\[
 H^1(Y_C)=q_C^*H^1(\Pic^d C),\quad
 H^2(Y_C)=q_C^*H^2(\Pic^d C)\oplus\Q\tau_C.
\]
The second summand may use $\tau$ because its coefficient on
$c_1(\det S_2)$ is the nonzero integer $u$.

\begin{proposition}\label{prop:pencil-pic}
For $g\ge6$, $d\le g-1$ and $\rho\ge3$,
\[
 H^1(\Pencils;\Q)=0,\qquad
 \Pic(\Pencils)_\Q\xrightarrow[\sim]{c_1}H^2(\Pencils;\Q)
                    =\Q\lambda\oplus\Q a_0\oplus\Q b_0.
\]
\end{proposition}
\begin{proof}
Pullback and the global class $\tau$ define a morphism of
constructible complexes on $\Mg$,
\[
       Ra_*\Q_{\Jac}\oplus\Q_{\Mg}[-2]
                           \longrightarrow Rp_*\Q_{\Pencils}.
\]
Proper base change and \cref{lem:fl} say it is an isomorphism on
cohomology sheaves in degrees zero, one and two. Their truncations
through degree two are therefore quasi-isomorphic. Taking
hypercohomology of these truncations gives
\[
 H^1(\Pencils)=H^1(\Jac),\qquad
 H^2(\Pencils)=H^2(\Jac)\oplus\Q\tau.
\]
The rigidified Jacobian calculation of Ebert--Randal-Williams
\cite[Theorem E and Theorem 4.6]{ERW} gives, for $g\ge6$,
$H^1(\Jac)=0$ and $H^2(\Jac)=\Q\lambda\oplus\Q\eta_0$.
These two classes have different geometric roles: $\lambda$ comes
from the source curve, while $\eta_0$ has a nonzero theta component on
the Jacobian fibre. The additional class $\tau$ records the choice of
a two-dimensional space of sections over that Jacobian.
The tangent Euler-class convention in that calculation changes the
sign of the $\beta$ term; its second generator is a nonzero rational
multiple of $-\eta_0$. For the algebraic Picard group and its comparison
with cohomology, see
\cite[Theorem B(i), Theorem 1.1(ii), and Corollary 1.2]{MV}.
Equation \eqref{eq:tau-eta} gives the displayed three Chern classes.

Apply \cref{lem:pic-c1} using the finite source-level presentation
from \cref{lem:full-pencil}. It supplies injectivity for the whole
rational Picard group; the three constructed lines supply surjectivity.
This uses vanishing of the invariant $H^1$ of the level cover, not
vanishing of its entire $H^1$.
\end{proof}

\subsection{Three relations on the simple locus}
We now turn the three generators into relations. Basepoint-free
evaluation supplies the first one. The other two come from
discriminants: one detects a repeated ramification point on the source,
and the other detects coincident branch values on the target. Both
discriminants are nowhere zero on the simple Hurwitz locus, so their
lines become trivial there.

The following relations hold on $\U_{d,g}$ whenever $g\ge2$, $d>0$
and the simple open is nonempty. They do not assume Picard generation.
Write $n=2g+2d-2$ in this geometric discussion.

\begin{lemma}\label{lem:relations}
In $\Pic(\U_{d,g})_\Q$,
\begin{equation}\label{eq:three-relations}
 a_0=0,\qquad 24\lambda+4a_0+6b_0=0,
                    \qquad2(n-1)(a_0+b_0)=0.
\end{equation}
\end{lemma}
\begin{proof}
Basepoint-free evaluation gives
\[
 0\longrightarrow\pi^*\det V\otimes L^{-1}
 \longrightarrow\pi^*V\longrightarrow L\longrightarrow0.
\]
Push its second Chern identity down the curve to obtain
$\alpha-d\gamma=0$.

The Wronskian is a section of
$\mathcal R=\omega_\pi L^2(\det V)^{-1}$, of class
$r_0=k+2l-\gamma$. Its zero scheme $R_f$ is finite flat of degree $n$:
on every fibre the section is nonzero and hence a nonzerodivisor.
The trace discriminant of the finite algebra $\pi_*\mathcal O_{R_f}$
is a section of $(\det\pi_*\mathcal O_{R_f})^{-2}$. GRR and
$\pi_*k^2=12\lambda$ give its class
\[
 \pi_*(r_0^2+r_0k)
 =24\lambda+4\alpha+6\beta-(6g-6+4d)\gamma
 =24\lambda+4a_0+6b_0.
\]
On the simple locus $R_f$ is \'etale over the base, so this section
is nowhere zero.

Finally use $P=\PP_{\rm lines}(V)$, with $z_0=c_1(\mathcal O_P(1))$.
On $\mathcal C\times P$ the universal linear combination is a section
of $L\mathcal O_P(1)$. Its zero scheme $Z_f\to P$ is finite flat
of degree $d$, even before excluding base points, since no nonzero
section of the pencil vanishes identically on a fibre. Its trace
discriminant line has class
\[
                       \alpha+\beta+n z_0.
\]
It is $\mathcal O_P(n)\otimes\mathcal L$ for a line $\mathcal L$
on the scalar gerbe with $c_1(\mathcal L)=\alpha+\beta$. The trace
determinant is a binary form in $\operatorname{Sym}^nV^*\otimes\mathcal L$.
Its discriminant is a section of
\[
                  \mathcal L^{2(n-1)}(\det V)^{-n(n-1)}.
\]
The negative determinant exponent follows by multiplying the squared
wedges of the covector factors of a binary form. The scalar weights
cancel, so this line descends, with class $2(n-1)(a_0+b_0)$.
On the basepoint-free locus the incidence cover is the original cover
over its target pencil; simple ramification gives simple discriminant
zeros, and distinct branch values make the binary form squarefree.
Its discriminant is therefore nowhere zero. This also handles descent
to targets that are conics.
\end{proof}

The relation matrix in columns $(\lambda,a_0,b_0)$ has determinant
$-48(n-1)\ne0$. Restriction from the smooth full pencil stack to its
simple open is surjective on rational Picard groups by
\cref{lem:pic-c1}. Thus \cref{prop:pencil-pic,lem:relations} prove
\begin{equation}\label{eq:geometric-a}
                 \Pic(\U_{d,g})_\Q=0
       \quad(d\ge6,\ d+1\le g\le2d-5).
\end{equation}
In particular this proves the two needed cases $(6,7)$ and $(7,8)$.

\section{The surface case: a K3 test family}
\label{sec:k3}

In this section $d\ge6$, $g=2d-4$, and $\rho=2$. We prove
\[
           \Pic(\Pencils)_\Q=\Q\lambda\oplus\Q a_0\oplus\Q b_0,
\]
and then apply the three relations already established. The required
case is $(d,g)=(6,8)$. For $\rho=2$, the comparison in
\cref{lem:fl} provides only an injection in degree two. We use the
Picard group of the whole K3 pencil family to control restrictions of
global line bundles to a general fibre.

The difficulty at this endpoint is that the pencil scheme on a general
curve is a surface. The degree-two Lefschetz comparison does not rule
out additional divisor classes on that surface. We therefore begin
with a line bundle on the full pencil stack and ask what its fibre
classes can be. A K3 family will determine those classes. Since they
come from a global line bundle, they vary as a flat cohomology class,
and the information obtained from a test fibre can be transported over
the Petri-general locus.

Two algebraic steps then remain. Vanishing of the generic Chern class
leaves a possible algebraically trivial line bundle, which is handled
using the strong Franchetta theorem. A generically trivial line is
represented by a vertical divisor, so we must also show that such
divisors come from the moduli of source curves. We establish this last
fact first, before constructing the test family.

\subsection{Conormal spaces and vertical divisors}
The differential of the map forgetting a pencil is governed by the
Petri map. For a pencil, the kernel of the Petri map has a concrete form:
the basepoint-free pencil trick expresses it as sections of a line
bundle, and the map to quadratic differentials is multiplication by a
Wronskian. This description will let us bound the critical locus by
fixing one nonzero quadratic differential on the source curve.

First allow general $g\ge2,d>0$ and write $m=3g-3$.
For a pencil $x=(C,L,V)$, let $E$ be its full fixed divisor and put
$M_0=L(-E)$. Write $s_i=b a_i$, where $b$ cuts out $E$ and
$a_0,a_1$ span the moving basepoint-free pencil in $M_0$.

\begin{lemma}\label{lem:conormal}
Up to the determinant factor of the fixed two-space $V$, the conormal
space of $p:\Pencils\to\Mg$ at $x$ is
\begin{equation}\label{eq:conormal}
 H^0(K_C M_0^{-2}(-E))\longrightarrow H^0(K_C^2),
 \qquad t\longmapsto b t\,(a_0\,da_1-a_1\,da_0).
\end{equation}
This map is injective, and its image is $\ker(dp_x^*)$.
\end{lemma}
\begin{proof}
Tensor the basepoint-free evaluation sequence of the moving pencil by
$K_C L^{-1}=K_C M_0^{-1}(-E)$. The pencil trick identifies the kernel
of $V\otimes H^0(K_C L^{-1})\to H^0(K_C)$ with the domain in
\eqref{eq:conormal}. The symbol sequence of $\mathcal D(L)$ and
the deformation complex from \cref{lem:full-pencil} identify the dual
cokernel of $dp$ with this Petri kernel. More explicitly, the map
$H^1(\mathcal D(L))\to V^*\otimes H^1(L)$ is surjective, because
its sheaf map is injective with torsion cokernel. Adjusting a lift of
a curve deformation by $H^1(\mathcal O_C)$ gives the Petri cokernel
as the obstruction to lifting it with its pencil.

The induced dual map to quadratic differentials is computed locally
on a kernel element by
\[
 (a_1t)d(ba_0)-(a_0t)d(ba_1)
                         =bt(a_1da_0-a_0da_1).
\]
The $db$ terms and changes of line-bundle trivialization cancel. Up
to sign this is \eqref{eq:conormal}. Both $b$ and the moving Wronskian
are nonzero, proving injectivity. Incomplete pencils are included.
\end{proof}

\begin{lemma}\label{lem:fixed-cotangent}
For a fixed curve $C$ and nonzero $q_0\in H^0(K_C^2)$, only finitely
many degree-$d$ pencils satisfy $q_0\in\ker(dp_x^*)$.
\end{lemma}
\begin{proof}
If $R$ is the divisor of the moving Wronskian, \cref{lem:conormal} gives
$R+E\le\divisor(q_0)$. There are finitely many such effective divisor
pairs. Also $M_0^2=\mathcal O_C(R)K_C^{-1}$, so there are finitely
many possible $M_0$, differing by two-torsion.

For fixed $M_0$, the Wronski morphism
\[
 \Gr(2,H^0(M_0))\longrightarrow\PP H^0(K_C M_0^2)
\]
is everywhere defined: the Wronskian of any independent pair is
nonzero in characteristic zero. It is a linear subsystem of the
Pl\"ucker system and pulls $\mathcal O(1)$ back to the ample Pl\"ucker
line. No positive-dimensional projective fibre can have this property.
Thus the morphism is finite onto its image, and a specified Wronskian
divisor gives finitely many pencils. Finally $L=M_0(E)$.
\end{proof}

\begin{proposition}\label{prop:vertical}
If $\rho\ge2$, every vertical prime divisor of $\Pencils$ is the
unique reduced pullback of a prime divisor of $\Mg$. No prime divisor
lies over a base locus of codimension at least two.
\end{proposition}
\begin{proof}
Work on a smooth \'etale chart of $\Mg$ of dimension $m$.
For an irreducible locally closed corank-$c$ stratum $Y$, $c\ge1$,
let $Z$ be the closure of its image. Generic smoothness gives
$dp(TY)=TZ$ on dense smooth opens. The projectivized rank-$c$
conormal bundle over $Y$ maps to $\PP N^*_{Z/\Mg}$, whose dimension
is $m-1$. Its fibres are finite by \cref{lem:fixed-cotangent}.
Therefore
\[
                   \dim Y+c-1\le m-1,
             \qquad\dim\Crit(p)\le m-1.
\]
Every component of the smooth total space, of dimension $m+\rho$,
dominates the base, since otherwise its generic point would be critical.

Over the DVR at the generic point of a base divisor, the total space
is torsion-free and hence flat: the base uniformizer is nonzero on
each dominating component. Its fibre is a pure $\rho$-dimensional
local complete intersection, hence $S_2$. The critical locus over that
generic point is finite, by the preceding dimension bound. Thus the
geometric fibre is regular in codimension one for $\rho\ge2$, and
is normal by Serre's criterion. Fulton--Lazarsfeld connectedness
\cite[Theorem I]{FL}, followed by the connected Grassmannian fibres
of $G^1_d(C)\to W^1_d(C)$, gives geometric connectedness. A connected
normal scheme is integral. In particular there is one prime divisor
over the base divisor, with multiplicity one.

If a prime divisor of the total space mapped into a locus of codimension
at least two, generic smoothness on that divisor would give
$\rk(dp)\le\dim(\text{image})+1\le m-1$. It would be contained
in the critical locus, although its dimension is $m+\rho-1>m-1$,
a contradiction. This also shows that the unique pullback has no
additional prime components over smaller loci. The assertions descend
from the chart. Normality, not connectedness alone, is the step that
rules out splitting over a base divisor.
\end{proof}

\subsection{Pencils on a fixed K3 surface}
Return to $g=2d-4$. Choose a polarized K3 surface $T$ with
\[
             \Pic(T)=\Z H,\qquad H^2=2g-2=4d-10.
\]
Such primitive polarized surfaces exist. Every member of $|H|$ is
reduced and irreducible, since $H$ generates the Picard group.
Lazarsfeld's theorem \cite[Theorem, p.~299]{Lazarsfeld} gives a nonempty open of smooth
Petri-general curves. Choose a nonempty open $B_0$ in it and delete a
hyperplane in $|H|=\PP^g$; thus $\mathcal O_{|H|}(1)|_{B_0}$ is trivial.
Let $S_{B_0}\to B_0$ be the full relative degree-$d$ pencil scheme.

For these curves the Brill--Noether numbers are
\[
 \rho(g,1,d)=2,\quad \rho(g,1,d-1)=0,\quad
 \rho(g,1,d-2)=-2,\quad \rho(g,2,d)=2-d<0.
\]
Petri generality and nonemptiness at nonnegative $\rho$ imply that
$S_{B_0}$ is a smooth proper family of connected surfaces, all
degree-$d$ pencils are complete, and every based pencil has exactly
one simple fixed point. The degree-$(d-1)$ pencils form a nonempty
finite reduced scheme on each curve and are complete and basepoint-free.
The based locus in a surface is one-dimensional, so the basepoint-free
open $S_{B_0}^{\rm bf}$ is nonempty. Here nonemptiness at nonnegative
$\rho$ is the Brill--Noether existence statement obtained from
\cite[Theorem II(a) and \S2]{FL}; Petri generality gives smoothness and
the expected dimensions. The negative values above exclude a second
base point and exclude incomplete degree-$d$ pencils.

The next construction replaces a pencil on a moving curve by a vector
bundle on the fixed K3 surface. The curve is recovered as the zero
divisor of the determinant of two sections of that bundle. This allows
us to use the Picard group of a projective sheaf moduli space, whose
rank is known, to bound the Picard group of the whole pencil family.

For a complete basepoint-free pencil $A_C$ on a smooth $C\in|H|$,
define its Lazarsfeld--Mukai bundle by
\begin{equation}\label{eq:lm}
 0\longrightarrow E_{C,A_C}^*\longrightarrow
 H^0(A_C)\otimes\mathcal O_T\longrightarrow i_*A_C\longrightarrow0.
\end{equation}
The kernel is locally free, as seen in a local equation of the smooth
Cartier divisor $C$. Duality gives
\[
 0\to H^0(A_C)^*\otimes\mathcal O_T\to E_{C,A_C}
                   \to i_*(K_C A_C^{-1})\to0.
\]
Completeness and $H^1(\mathcal O_T)=0$ give
\begin{equation}\label{eq:lm-invariants}
 \begin{aligned}
 \det E&=H,\qquad c_2(E)=\deg A_C,\qquad h^0(E)=g+3-\deg A_C,\\
 H^0(E^*)&=H^1(E)=H^2(E)=0.
 \end{aligned}
\end{equation}
No assumption that the residual series is basepoint-free is needed.

Every such $E$ is slope-stable. A destabilizing rank-one subsheaf
would extend to a map $\mathcal O_T(kH)\to E$ with $k\ge1$.
Multiplication by a nonzero section of $(k-1)H$ would give a nonzero
section of $E(-H)=E^*$, contradicting \eqref{eq:lm-invariants}.
The rank-one slopes are integral multiples of $H^2$, whereas the
rank-two slope is $H^2/2$, so there is no equality wall. The same
observation makes every Gieseker-semistable sheaf of rank two and
first Chern class $H$ slope-stable.

Any two independent sections of a stable such bundle are generically
independent. Otherwise their saturated rank-one image extends to
$\mathcal O_T(kH)$; having two sections forces $k\ge1$, contradicting
stability. Their determinant is therefore never identically zero.
If its zero curve is smooth, the local section matrix has rank one
along that curve and reduces to $\operatorname{diag}(1,t)$.
Its cokernel is a line on the curve. Dualizing reconstructs
\eqref{eq:lm}, and $H^0(E^*)=H^1(E^*)=0$ recovers the complete
two-space of sections. This is an inverse in families: the cokernel
is a relative line bundle, and the operations commute with base change.
Simultaneous scalar twists of a universal bundle and its section
space do not change the pencil. Thus the inverse works also when
only a twisted universal sheaf exists.

\subsection{Two rank drops in the sheaf-moduli construction}
There will be two restrictions to open sets. First we discard the
non-locally-free sheaves; their divisor removes a nonzero Picard class.
After adding a two-space of sections, we restrict to curves in $B_0$.
The removed hyperplane in $|H|$ makes the determinant line trivial and
removes another class. We keep track of both restrictions because the
intermediate Grassmannian bundle itself contributes one new class.

On a K3 surface the Mukai vector of a bundle is
\[
 v(E)=\bigl(\rk E,c_1(E),\rk E+\tfrac12c_1(E)^2-c_2(E)\bigr).
\]
The Mukai pairing is
$(r,c,s)\cdot(r',c',s')=c\cdot c'-rs'-r's$.
For the degree-$d$ bundles this gives
\[
                   v=(2,H,d-3),\qquad v^2=2.
\]
It is primitive, and $H$ is generic for it by the absence of slope
equalities. The moduli space $\mathcal N=M_H(v)$ is a smooth
projective irreducible symplectic fourfold. The Mukai Hodge comparison
\cite[Theorem 0.1]{Yoshioka} identifies its $H^2$ with $v^\perp$ in
the Mukai lattice $H^0(T,\Z)\oplus H^2(T,\Z)\oplus H^4(T,\Z)$.
The algebraic Mukai lattice of $T$ has rank three, and its intersection
with $v^\perp$ has rank two. Since $\Pic^0(\mathcal N)=0$,
\begin{equation}\label{eq:moduli-pic}
                         \dim_\Q\Pic(\mathcal N)_\Q=2.
\end{equation}

There is a unique stable sheaf with vector
$v_0=(2,H,d-2)$, of square $-2$, and it is a locally free bundle $E_0$.
Existence follows from a degree-$(d-1)$ pencil on one of the chosen
curves. For uniqueness, two nonisomorphic stable sheaves with this
vector have no Hom in either direction and hence no $\Ext^2$; but
their Euler pairing is $-v_0^2=2$, a contradiction. Local freeness
also follows by taking the double dual: a singularity of length
$l>0$ would give a stable double dual of square $-2-4l$, contradicting
$\dim\Ext^1=v^2+2\ge0$. By \eqref{eq:lm-invariants}, $h^0(E_0)=d$.

Now let $F\in\mathcal N$ be non-locally-free, with
$l=\operatorname{length}(F^{**}/F)\ge1$. The double dual remains
stable, since a destabilizing line subsheaf would intersect $F$ in
one of the same slope. Its vector is $v+(0,0,l)$ and its square is
$2-4l\ge-2$. Thus $l=1$ and $F^{**}=E_0$. Conversely every kernel
of $E_0\to k(x)$ is stable of vector $v$. The point and quotient
line are recovered from the kernel up to scalar automorphisms of
$E_0$. Hence the nonfree locus is the image of the integral
three-dimensional projective bundle of quotient lines of $E_0$ over
$T$, with finite fibres. Its reduced support is one irreducible
divisor $D_{\rm nf}$. It is nonzero in $\Pic(\mathcal N)_\Q$, because
it has positive intersection with the cube of an ample class.

Let $T_0\subset\mathcal N$ be the nonempty open of locally free
bundles with $H^1(E)=H^2(E)=0$. Divisor localization and
\eqref{eq:moduli-pic} show
\begin{equation}\label{eq:first-drop}
                              \dim_\Q\Pic(T_0)_\Q\le1.
\end{equation}
Form the possibly twisted Grassmannian $P\to T_0$ of two-spaces in
$H^0(E)$; its fibres are $\Gr(2,d-1)$. A line whose restrictions to
all fibres are trivial is pulled back from $T_0$, by cohomology and
base change and the evaluation map. Fibre degrees contribute at most
one additional generator. Possible Brauer obstructions only restrict
the integral fibre degrees, so
\[
                              \dim_\Q\Pic(P)_\Q\le2.
\]
The determinant of two sections defines an everywhere-defined morphism
$f:P\to|H|$. Its line $\delta=f^*\mathcal O(1)$ restricts to
Pl\"ucker on every Grassmannian fibre, and is nonzero rationally.
The family-level inverse above identifies
$f^{-1}(B_0)=S_{B_0}^{\rm bf}$. On this open $\delta$ is trivial.
Every line on the open extends to the smooth $P$, giving the second
rank drop
\begin{equation}\label{eq:second-drop}
                        \dim_\Q\Pic(S_{B_0}^{\rm bf})_\Q\le1.
\end{equation}
The Grassmannian fibres used in this calculation belong to $P\to T_0$.
They show that $\delta$ is nonzero before restriction, whereas the
choice of $B_0$ shows that it is zero afterwards. These are the two
facts needed for the second rank drop.

\subsection{Restoring the one basepoint divisor}
So far we have bounded the Picard group of the basepoint-free family.
To return to all pencils we must count the divisorial components that
were omitted. Every based pencil here is obtained by adding one point
to a degree-$(d-1)$ pencil. The bundle $E_0$ will show that this locus
has just one irreducible component, so restoring it adds at most one
Picard generator.

The determinant construction for $E_0$ is a morphism
\[
                      f_0:\Gr(2,d)\longrightarrow|H|=\PP^g.
\]
It pulls back $\mathcal O(1)$ to the ample Pl\"ucker line, so it is
finite onto its image. Both spaces have dimension $2(d-2)=g$,
and therefore $f_0$ is surjective. Over $B_0$ its inverse construction
parametrizes all degree-$(d-1)$ pencils, since their stable bundles
are all $E_0$. The inverse image $f_0^{-1}(B_0)$ is irreducible.
Pulling back the universal smooth curve and adding its chosen point
to the minimal pencil gives precisely the entire based locus in
$S_{B_0}$. There is only one simple basepoint and no lower-degree
pencil. Thus this locus has one irreducible divisorial support.
Restoring it in \eqref{eq:second-drop} gives
\[
                              \dim_\Q\Pic(S_{B_0})_\Q\le2.
\]

These two dimensions are attained by $a_0,b_0$. On a fixed Petri fibre,
normalize the Poincar\'e bundle at a point. Then
$\alpha=-2q^*\theta$ and $\beta=0$: the square of the mixed K\"unneth
component integrates to $-2\theta$, whereas multiplication by $k$
has no surviving mixed term. The $\rho=2$ version of \cref{lem:fl}
injects the ambient $H^2$ into that of the pencil surface, so
$q^*\theta$ and $\gamma=c_1(\det S_2)$ remain independent. The
restrictions
\[
                     a_0=-2q^*\theta-d\gamma,
               \qquad b_0=-(g-1)\gamma
\]
are therefore independent. We have proved the family calculation
\begin{equation}\label{eq:k3-family-pic}
                        \Pic(S_{B_0})_\Q=\Q a_0\oplus\Q b_0.
\end{equation}

\subsection{Global line bundles and the generic \texorpdfstring{$\Pic^0$}{Picard variety}}
We now apply the family calculation to a line bundle that is already
defined on $\Pencils$. The K3 family supplies its two fibre coefficients.
The remainder will first be shown to be trivial at the generic point,
and only then expressed as a pullback from $\Mg$.

Let $\mathcal L\in\Pic(\Pencils)_\Q$. Subtract the unique rational
combination of $a_0,b_0$ making its restriction to $S_{B_0}$ zero,
and call the result $\mathcal L'$. On the connected Petri open
$Y\subset\Mg$, the map $p$ is smooth proper. Fibre Chern classes
of a global line form a flat section of $R^2p_*\Q$. Since the section
vanishes on one test fibre, it vanishes throughout $Y$. It follows,
after a tensor power, that the geometric generic restriction of
$\mathcal L'$ lies in $\Pic^0$. This can equally be formulated in
the relative Picard scheme: zero rational fibre Chern class means
membership in $\Pic^\tau$, whose quotient by $\Pic^0$ is torsion.
No dominance of the K3 family over $\Mg$ is needed; the line being
tested was global before restriction.

Let $F=\C(\mathcal M_g)$ and $C_\eta$ be its generic curve. The
$H^1$ isomorphism of \cref{lem:fl} gives an $F$-defined isogeny
\[
       \Pic^0(\Pic^d(C_\eta))\longrightarrow\Pic^0(G^1_d(C_\eta)).
\]
The source is the dual Jacobian, canonically a Jacobian via its
principal polarization, even though $\Pic^d(C_\eta)$ is a torsor.
The strong Franchetta theorem \cite[Theorem 5.1]{Schroer}, for the
generic curve of genus $g\ge3$ with no marked points, says that the
$F$-rational points of its Picard scheme are generated by the canonical
class. Its degree-zero part therefore gives
$J(C_\eta)(F)\otimes\Q=0$. An isogeny has an $F$-defined inverse
up to multiplication by an integer; hence
\[
                    \Pic^0(G^1_d(C_\eta))(F)\otimes\Q=0.
\]
Actual line bundles inject into the $F$-points of the Picard scheme;
the Brauer obstruction lies on the cokernel side. Thus
$\mathcal L'_\eta=0$ in the rational Picard group. The geometric
abelian variety $\Pic^0$ has not been set to zero.

A line generically trivial on the regular total stack is represented
rationally by a vertical divisor, by extending its generic
trivialization. \Cref{prop:vertical} identifies every such prime
with a reduced base pullback. The base Picard group is
$\Pic(\Mg)_\Q=\Q\lambda$. One can obtain this algebraic statement
from \cite[Theorem B(i) and Corollary 1.2]{MV}: pullback to the proper connected universal Jacobian
is injective, and its second generator has a nonzero theta restriction,
so only the Hodge direction can come from the base. We conclude
\[
                     \mathcal L\in\Q\lambda+\Q a_0+\Q b_0.
\]
The two fibre classes are independent, and properness with geometrically
connected fibres gives $p_*\mathcal O=\mathcal O$, so pullback of
the nonzero Hodge class is injective. This proves the asserted direct
sum for $\Pic(\Pencils)_\Q$.

Finally restrict to $\U_{d,g}$ and use \cref{lem:relations}. Its
nonzero relation determinant kills all three generators. Thus
\begin{equation}\label{eq:rho-two-result}
                        \Pic(\U_{d,2d-4})_\Q=0\qquad(d\ge6).
\end{equation}
Together with \eqref{eq:geometric-a}, this proves
\cref{thm:geometric}.

\section{Completion of the proof}
\label{sec:completion}

\begin{proof}[Proof of \cref{thm:main}]
Suppose first that $n=2g+2d-2\ge3$.
For $g\ge3$ and $2\le d\le5$, apply \cite[Theorem A]{DP}, with
its unparameterized-target and unordered-branch convention.
For $d\ge g$, apply the unmarked-branch assertion of
\cite[Theorem 1.1]{Mullane}. This includes the low genera $g=1,2$;
there is no use of the different inequality for its one-marked-branch
variant.

For completeness, genus zero can also be handled directly. The space
of degree-$d$ rational functions is the open in $\PP^{2d+1}$ of pairs
of binary degree-$d$ forms with nonzero resultant. Restrict further
to simple covers. Divisor localization shows that its rational Picard
group vanishes, since the nonzero resultant divisor kills the
hyperplane generator. Source and target coordinate changes are by
$\PGL_2\times\PGL_2$, a connected group with no rational characters.
The equivariant Picard sequence therefore gives the same vanishing
for its quotient. This is the unparameterized genus-zero Hurwitz stack.

It remains to consider $d\ge6$ and $g>d$. If $g\ge9$, use
\cref{thm:uniform}. If $g\le8$, the only possibilities are
$(6,7),(6,8),(7,8)$, all covered by \cref{thm:geometric}.
Thus the parameter ranges are exhaustive.

Finally $n<3$ with $d\ge2$ forces $(d,g)=(2,0)$. Every such map is
isomorphic to $z\mapsto z^2$. Its automorphism group, allowing both
source and target automorphisms, is $\Gm\rtimes\mu_2$, with the
involution acting on $\Gm$ by inversion. Its character group has
zero rationalization, so its classifying stack has rational Picard
group zero; its coarse space is a point. This treats the non-finite
inertia case separately. For $n\ge3$, inertia is finite, and the
rational coarse-space assertion follows by descent as above.
\end{proof}

{\raggedright
\bibliographystyle{alpha}
\bibliography{references}
}
\par\bigskip
\begingroup
\footnotesize
\noindent University of Science and Technology of China\par
\noindent\textit{Email: }\href{mailto:wenhailin@mail.ustc.edu.cn}{\texttt{wenhailin@mail.ustc.edu.cn}}\par
\endgroup

\end{document}